\pdfoutput=1
\documentclass[11pt]{article}

\usepackage{amsmath}
\usepackage{amssymb}
\usepackage{amsfonts}
\usepackage[table]{xcolor}
\usepackage[bottom]{footmisc}
\usepackage{graphicx}
\usepackage[hyphens]{url}
\usepackage{multirow}
\usepackage{array}
\usepackage{tabularx}
\usepackage{booktabs}

\DeclareMathOperator*{\Minimize}{Minimize}
\DeclareMathOperator*{\argmin}{argmin}

\usepackage{fontawesome5}
\usepackage{pifont}
\newcommand{\tY}{\textcolor[rgb]{0,0.55,0}{\ding{51}}}
\newcommand{\tN}{\textcolor[rgb]{0.8,0,0}{\ding{55}}}
\newcommand{\tP}{\textcolor[rgb]{0.9,0.55,0}{\faAdjust}}
\usepackage{siunitx}

\usepackage[english]{babel}
\usepackage{csquotes}
\usepackage[style=apa,backend=biber]{biblatex}
\AtBeginBibliography{\setlength{\emergencystretch}{1.5em}}
\usepackage[hidelinks]{hyperref}
\begin{document}

\title{Solving the capacitated vehicle routing problem with polygonal neighborhoods using a metaheuristic enhanced with block coordinate descent}

\author{
    Ernesto G. Birgin\thanks{Department of Computer Science, Institute of
    Mathematics and Statistics, University of S\~ao Paulo, Rua do
    Mat\~ao, 1010, Cidade Universit\'aria, 05508-090, S\~ao Paulo, SP,
    Brazil. e-mail: egbirgin@ime.usp.br. ORCID: 0000-0002-7466-7663. Corresponding author.}
\and 
    Rogério S. Matos\thanks{Department of Applied Mathematics, Institute of
    Mathematics and Statistics, University of S\~ao Paulo, Rua do
    Mat\~ao, 1010, Cidade Universit\'aria, 05508-090, S\~ao Paulo, SP,
    Brazil. e-mail: rogerio.matos@ime.usp.br}
}

\date{June 29, 2026\thanks{Revision made on November 11, 2026.}}

\maketitle

\begin{abstract}

We study a capacitated vehicle routing problem (CVRP) in which customers are represented by planar regions, called customer neighborhoods, rather than fixed points. Visiting a customer consists of selecting a point within its neighborhood. Therefore, besides determining the sequence of customers visited by each vehicle, one must also determine the corresponding visit points. To better reflect practical applications, customer neighborhoods are allowed to be simple polygons, either convex or non-convex, and distances are measured using the Euclidean norm. Applications include aerial reconnaissance, wireless sensor data collection and inspection, and two-stage delivery systems. We propose a general variable neighborhood search framework in which each neighbor is represented by a collection of routes. Constructing a route involves two stages, defining a customer sequence and then determining the associated visit points. The latter task is solved by a block coordinate descent (BCD) procedure, the central methodological contribution of this work, which computes high-quality visit points at a cost low enough to be used repeatedly inside the search. In extensive experiments on newly generated and adapted benchmark instances, within the time limits considered, the proposed framework matches all solutions proven optimal by CPLEX and matches or improves the best solutions found by CPLEX in the remaining instances. For the instances adapted from the CVRP literature, it also improves, in every instance, a CVRP-based upper bound given by the solution of the original CVRP instance reported in the literature, which is feasible for the problem with neighborhoods.\\

\end{abstract}

\section{Introduction}

In the classic vehicle routing problem \parencites[VRP;][]{Dantzig1959}[see also][]{Golden2023,Laporte2009,Toth2014}, there are $m$ customers, a depot, and an unlimited number of identical trucks. The problem consists in building a set of routes such that each route starts and ends at the depot and the routes, as a whole, visit each customer exactly once. The combined length of the set of routes must be minimized. Note that in the classical problem customers have fixed geographical locations such that the distance between each pair of customers or between a customer and the depot is known in advance and fixed. Thus, the length of a route that leaves the depot, visits customers $i_1, i_2, \dots, i_{\kappa}$, and returns to the depot is given by the simple sum of the distance from the depot to customer $i_1$, the distances between consecutive customers along the route, and the distance between customer $i_{\kappa}$ and the depot.

In the VRP with neighborhoods (VRPN), each customer is represented by an arbitrary region of the plane, termed customer neighborhood, rather than by a coordinate location point. This idea was introduced in the context of the traveling salesman problem (TSP) by~\textcite{Arkin1994}. In this context, visiting a customer means visiting a point within the customer's neighborhood. In addition to the challenges of the VRP, the visit point for each customer must be determined. Furthermore, because the location of visit points is a variable in the problem, the distances between pairs of customers and between the customers and the depot are not fixed. In the specific case when customer neighborhoods are given by two-dimensional balls (circles), the problems are known as Close Enough TSP \parencite[CETSP;][]{Behdani2014,BoYuan2007,Carrabs2017,Carrabs2020,Coutinho2016,DiPlacido2022,DiPlacido2023,Gulczynski2006,Lei2024,Wang2019,Zhang2023} or Close Enough VRP \parencite[CEVRP;][]{Mennell2009,Procopio2020}. The names come from the fact that the radius of a customer's circle is interpreted as the maximum Euclidean distance from the customer at which the route can pass for it to be considered as visited. Neighborhoods other than disks have also been considered for the TSP and related problems. \textcite{Gentilini2013} proposed an MINLP formulation and an exact method for polyhedral and ellipsoidal neighborhoods, \textcite{Puerto2022} proposed conic formulations and a heuristic for routing drones through convex sets and polygonal chains, and \textcite{Kulich2022} considered polygonal neighborhoods among polygonal obstacles. For the VRP, neighborhoods given by balls were considered by~\textcite{Mennell2009,Procopio2020}, with applications in aerial reconnaissance route planning. Because disks are convex, in these works the visit points of a fixed sequence of customers can be computed exactly by solving a second-order cone program for the whole route~\parencite[see, for example,][]{Coutinho2016,Mennell2009}. Polygonal neighborhoods, possibly non-convex, were considered together with capacity constraints by~\textcite{Gao2021}, under the name capacitated vehicle routing with target geometric constraints (CVRG), in a robotic delivery setting, where the visit points of a fixed sequence are improved by moving one point at a time. They proposed an exact dynamic programming algorithm, whose optimality is guaranteed for convex neighborhoods, and a heuristic based on dynamic programming.

\begin{figure}[ht!]
\centering
\includegraphics[angle=90,scale=1.0]{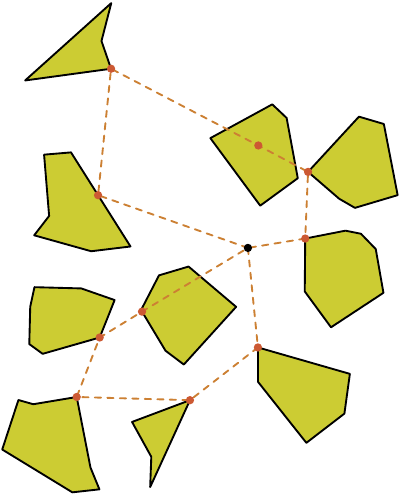}
\caption{Example of a two-route optimal solution to a ten-customer instance of the capacitated vehicle routing problem with non-convex polygonal customer neighborhoods. Note that each customer's visit point may be at a vertex, on an edge, or even within the interior of the polygon representing the customer's neighborhood. Because these locations directly dictate route length, the distance for a given visit sequence can only be calculated once the visit points are selected.}
\label{fig1}
\end{figure}

In the present work, we consider a VRPN in which each customer's neighborhood is represented by a simple polygon, possibly non-convex. See Figure~\ref{fig1} for an illustrative example, in which the visit points lie at vertices, on edges, or in the interior of the polygons. We consider that each customer $i$ has a demand $d_i>0$ and that all trucks have the same maximum capacity $\mathcal{C}$. Therefore, we call it the capacitated VRPN (CVRPN) problem. We also consider the Euclidean norm as a measure of distance; this implies that the problem is symmetric. The CVRPN that we study here includes the practical case in which trucks first deliver the aggregate demand of a neighborhood to a drop-off point located anywhere within that neighborhood. Crowdsourced couriers then complete the last mile~\parencite{Alnaggar2021,Le2019}. Each courier operates within a specific region represented by a non-convex polygon reflecting realistic urban boundaries or accessibility constraints. For instance, a courier's region may be a set of blocks bounded by a river, a railway, or an expressway, or a residential complex with an L- or U-shaped layout. The problem under consideration also encompasses scenarios in which vehicles are required to visit customers, not to deliver goods, but rather to perform measurements, such as collecting or inspecting data via wireless sensors. This technology, called automated meter reading \parencite[AMR;][]{Dong2007,Shuttleworth2008}, is used by utility service providers, such as electricity, water, and gas. It enables remote monitoring, billing, and management. In this scenario, it is unnecessary for the vehicle to reach the customer's exact location; it suffices for the vehicle to pass close enough for the measurement to be completed. Although such proximity is often modeled as a circular region around the customer, signal propagation may be obstructed in practice by walls, buildings, or other obstacles~\parencite{Haslett2007}, resulting in irregular coverage zones. For instance, the region from which a meter at the rear of a building can be read may consist of two separate street segments connected only through a side passage. Consequently, the effective measurement region is more accurately represented by a non-convex polygon, which captures the true spatial constraints imposed by the environment on the quality of the signal transmission. In all these cases, approximating the service region by a circle or a convex polygon would either include points that cannot serve the customer or exclude points that can, whereas arbitrary polygons allow the region to be represented as it is. Finally, the CEVRP in which the distances that define the neighborhoods are measured using the Manhattan or Chebyshev norms is a particular case of the CVRPN with convex polygonal customer neighborhoods, because these neighborhoods are then squares.

We say that a sequence that starts at the depot, is followed by a sequence of customers $i_1, i_2, \dots, i_{\kappa}$ and ends at the depot is a \textit{visit sequence}. In order for a visit sequence to become a route, we need to find the visit points for each of the customers $i_1, i_2, \dots, i_{\kappa}$. In order for the route to have minimal length, the visit points must be a solution to the problem
\begin{equation} \label{proba}
\Minimize_{\{ p_{i_j} \in \mathbb{R}^2, j=1,\dots,\kappa \}} \sum_{j=0}^\kappa \| p_{i_{j+1}} - p_{i_j} \| \mbox{ subject to } p_{i_j} \in P_{i_j} \mbox{ for } j=1,\dots,\kappa.
\end{equation}
In~\eqref{proba}, $p_{i_0}$ and $p_{i_{\kappa+1}}$ both represent the given geographical position $p_0$ of the depot, $P_{i_j}$ is the polygon representing the neighborhood of customer $i_j$, and the variable $p_{i_j}$ is the visit point, to be determined, of customer $i_j$. If the polygons representing the customer neighborhoods were all convex, problem~\eqref{proba} would be a nonlinear optimization problem. When non-convex polygons exist, problem~\eqref{proba} can be reformulated as a problem with a linear objective function, a convex feasible region determined by linear and quadratic constraints, and continuous and binary variables. This kind of problem is known as mixed-integer quadratically constrained programming (MIQCP). Commercial exact solvers such as CPLEX and Gurobi offer features to search for the global solution to this class of problems. In this paper, we propose solving problem~\eqref{proba} with a block coordinate descent method \parencite{Amaral2022,Birgin2022,Wright2015}. 

The flexibility gained by replacing fixed customer locations with polygonal neighborhoods comes at a significant computational cost. For a given set of routes, the corresponding visit points must also be determined, and these visit-point decisions directly affect the route lengths. Consequently, evaluating a candidate solution requires determining a set of customer sequences and, in addition, solving an MIQCP for each sequence. Because modern local-search and neighborhood-search methods evaluate a large number of candidate solutions, the efficient computation of visit points becomes a critical component of any practical algorithm for the CVRPN.

The CVRPN is NP-hard because it has the VRP as a special case. Therefore it is natural to address the CVRPN with a heuristic approach. In this paper, we consider a general variable neighborhood search (GVNS) method~\parencite{Mladenovi1997}, where, given a current solution, different search neighborhoods are inspected. Each search neighborhood is formed by neighbors, and each neighbor corresponds to a set of routes. Now, building each route by determining its visit points so its length can be computed involves solving a problem of the form~\eqref{proba}. Therefore, because we need to solve problem~\eqref{proba} repeatedly, solving~\eqref{proba} efficiently is crucial for the overall approach to be practical.

To the best of our knowledge, no previous work on the capacitated VRP with simple, possibly non-convex, polygonal neighborhoods provides a mathematical programming formulation and a procedure for computing the visit points in which the subproblem of each customer is solved exactly. The main contributions of this work are as follows. First, we study the capacitated vehicle routing problem with polygonal customer neighborhoods (CVRPN), a vehicle routing problem in which customers are represented by simple polygons, possibly non-convex. Second, we propose a mixed-integer quadratically constrained programming formulation for the problem. Third, we develop a block coordinate descent (BCD) procedure that efficiently computes visit points associated with a fixed sequence of customers. Fourth, we embed this route-evaluation procedure within a general variable neighborhood search (GVNS) metaheuristic specifically tailored to the CVRPN. Finally, we introduce benchmark instances and conduct extensive computational experiments that demonstrate the effectiveness of the proposed approach.

Table~\ref{tablit} compares the present work with representative works on routing problems with neighborhoods in terms of the presence of capacity constraints, the neighborhoods considered, the solution methods proposed, and the numerical experiments reported. In the table, TSPN denotes the TSP with neighborhoods, XPPN the crossing postman problem with neighborhoods studied by~\textcite{Puerto2022}, GCETSP the generalized CETSP, CVRG the problem studied by~\textcite{Gao2021}, and CVRPN the problem studied in this work. All works consider Euclidean distances in the plane. A check mark indicates that a work has the feature, a cross that it does not, and a half circle that the feature is present only partially or in restricted cases. Column ``Formulation'' indicates a mathematical programming formulation of the whole problem. Columns ``Visit points (solver)'' and ``Visit points (algorithm)'' indicate that the visit points of a fixed sequence are computed with an off-the-shelf solver, for example, by solving a second-order cone program with CPLEX or Gurobi, or with a specific algorithm described in the work, respectively. Column ``Exact'' indicates a method that guarantees an optimal solution of the whole problem, which in the present work corresponds to solving the MIQCP model of Section~\ref{sec2} with CPLEX, and column ``Heuristic'' a heuristic or approximation method. In the numerical experiments, ``Non-convex'' indicates experiments with non-convex neighborhoods, ``Benchmark'' publicly available benchmark instances, ``Literature'' a comparison with methods or results from the literature, and ``Extensive'' experiments on at least 200 instances, including instances with at least 100 customers. Half circles correspond, for example, to polygons handled only through their boundaries \parencite{Kulich2022,Yang2018}, non-convex polygons without a described method \parencite{Gao2021}, a formulation only of the subproblem with a fixed sequence \parencite{Coutinho2016,Lei2024,Zhang2023}, or experiments on at least 200 instances, all of them with fewer than 100 customers \parencite{Procopio2020,Puerto2022}. The present work is the only one in Table~\ref{tablit} that handles non-convex polygonal neighborhoods as regions, and it combines them with capacity constraints, a mathematical programming formulation, and a specific algorithm for computing the visit points.

\begin{table}[!tbp]
\centering
\footnotesize
\setlength{\tabcolsep}{2pt}
\begin{tabularx}{\textwidth}{>{\raggedright\arraybackslash}Xl*{13}{>{\centering\arraybackslash}p{0.62cm}}}
\hline
 & & & \multicolumn{3}{c}{Neighborhoods} & \multicolumn{5}{c}{Solution methods} & \multicolumn{4}{c}{Experiments} \\
\cmidrule(lr){4-6}\cmidrule(lr){7-11}\cmidrule(lr){12-15}
Works & Problem & \rotatebox{90}{Capacity constraints} & \rotatebox{90}{Disks} & \rotatebox{90}{Convex polygons} & \rotatebox{90}{Non-convex polygons} & \rotatebox{90}{Formulation} & \rotatebox{90}{Visit points (solver)} & \rotatebox{90}{Visit points (algorithm)} & \rotatebox{90}{Exact} & \rotatebox{90}{Heuristic} & \rotatebox{90}{Non-convex} & \rotatebox{90}{Benchmark} & \rotatebox{90}{Literature} & \rotatebox{90}{Extensive} \\
\hline
\nptextcite{Arkin1994} & TSPN & \tN & \tP & \tP & \tP & \tN & \tN & \tN & \tN & \tY & \tN & \tN & \tN & \tN \\
\nptextcite{Gentilini2013} & TSPN & \tN & \tY & \tY & \tN & \tY & \tY & \tN & \tY & \tY & \tN & \tY & \tP & \tN \\
\nptextcite{Yang2018} & TSPN & \tN & \tY & \tP & \tP & \tN & \tN & \tP & \tN & \tY & \tP & \tY & \tY & \tN \\
\nptextcite{Kulich2022} & TSPN & \tN & \tN & \tP & \tP & \tN & \tN & \tP & \tN & \tY & \tP & \tY & \tP & \tN \\
\nptextcite{Puerto2022} & XPPN & \tN & \tY & \tY & \tP & \tY & \tY & \tN & \tY & \tY & \tP & \tP & \tN & \tP \\
\nptextcite{Coutinho2016} & CETSP & \tN & \tY & \tN & \tN & \tP & \tY & \tN & \tY & \tN & \tN & \tY & \tY & \tY \\
\nptextcite{DiPlacido2023} & GCETSP & \tN & \tY & \tN & \tN & \tN & \tY & \tY & \tN & \tY & \tN & \tP & \tN & \tN \\
\nptextcite{Zhang2023} & CETSP & \tN & \tY & \tN & \tN & \tP & \tY & \tN & \tY & \tP & \tN & \tY & \tY & \tY \\
\nptextcite{Lei2024} & CETSP & \tN & \tY & \tN & \tN & \tP & \tY & \tY & \tN & \tY & \tN & \tY & \tY & \tN \\
\nptextcite{Mennell2009} & CEVRP & \tY & \tY & \tN & \tN & \tY & \tY & \tP & \tN & \tY & \tN & \tY & \tY & \tN \\
\nptextcite{Procopio2020} & CEVRP & \tP & \tY & \tN & \tN & \tY & \tN & \tP & \tY & \tY & \tN & \tY & \tP & \tP \\
\nptextcite{Gao2021} & CVRG & \tY & \tN & \tY & \tP & \tN & \tN & \tP & \tP & \tY & \tP & \tN & \tP & \tN \\
This work & CVRPN & \tY & \tN & \tY & \tY & \tY & \tY & \tY & \tY & \tY & \tY & \tY & \tY & \tY \\
\hline
\end{tabularx}
\caption{Representative works on routing problems with neighborhoods.}
\label{tablit}
\end{table}

The rest of this paper is organized as follows. Section~\ref{sec2} describes the problem and presents an MIQCP model. Section~\ref{sec3} explains how to transform a visit sequence into a route. This is accomplished by computing the customer visit points and determining the route's length using block coordinate descent. Section~\ref{sec4} describes the GVNS method and the search neighborhoods considered in this work. Section~\ref{sec5} presents exhaustive numerical experiments. Section~\ref{sec6} is devoted to conclusions and future work.

\vspace{2mm}
\noindent
\textbf{Notation.} We denote by $\mathbb{R}_+$ the set of non-negative real numbers. For a point $p$ in the plane, we denote its abscissa and ordinate by $p^x$ and $p^y$, respectively. The symbol $\| \cdot \|$ represents the Euclidean norm. Given $(x,y) \in \mathbb{R}^2$, we denote $(x,y)^{\perp} = (y,-x)$, that is, $(x,y)^{\perp}$ corresponds to a 90-degree clockwise rotation of $(x,y)$. Given two points $p$ and $q$ on the plane, the segment joining them is denoted by $[p, q]$.

\section{Problem description and model} \label{sec2}

The data for an instance of the CVRPN are: the position $p_0 \in \mathbb{R}^2$ of the depot, the maximum capacity $\mathcal{C} \in \mathbb{R}_+$, which is identical for all trucks in the fleet, the number of customers $m$ and, for each customer $i \in \{1,\dots,m\}$, its demand $d_i$, with $0 < d_i \leq \mathcal{C}$, and its polygonal neighborhood $P_i$. In this work, we consider polygonal neighborhoods given by polygons with at least three vertices. Without loss of generality, we assume that each polygon $P_i$ is given by the union of $\nu_i$ disjoint convex polygons, i.e., $P_i = \cup_{\ell=1}^{\nu_i} P_{i\ell}$, where the $P_{i\ell}$ are all convex. In other words, the collection $\{P_{i\ell}\}_{\ell=1}^{\nu_i}$ makes up a partition of $P_i$ into convex polygons. This partition is not unique. Any partition into convex polygons yields an equivalent model, but the number of binary variables $w_{i\ell}$ grows with $\nu_i$, so partitions with few polygons are preferable. Furthermore, we assume\footnote{This is the assumption that requires polygons to have at least three vertices. It would be simple to include polygons with one and two vertices as well, but they would have to be represented differently. If a polygon $Q$ is formed by a single vertex $v$, the constraint $q \in Q$ is written as $q=v$. If a polygon $Q$ has exactly two vertices $v_1$ and $v_2$, the constraint $q \in Q$ can be written as $q = v_1 + \lambda (v_2-v_1)$ plus $0 \leq \lambda \leq 1$. These details could easily be added to the model that follows, but we prefer to ignore these cases to simplify the presentation.} that each convex polygon $P_{i\ell}$ is expressed as $P_{i\ell} = \{ p \in \mathbb{R}^2 \;|\; a_{i\ell k} \, p^x + b_{i\ell k} \, p^y + c_{i\ell k} \leq 0 \mbox{ for } k=1,\dots,n_{i\ell} \}$, where $n_{i\ell}$ is the number of vertices of the polygon $P_{i\ell}$.

The variables of the model are summarized in Table~\ref{tabvar}. In all these variables, index~0 refers to the depot, whose position $p_0$ is fixed.

\begin{table}[ht!]
\centering
\begin{tabular}{lll}
\hline
Variable & Type & Main role \\
\hline
$x_{ij}$ & binary & equals 1 if customer $i$ immediately precedes customer $j$ on a route \\
$p_i$ & continuous & visit point of customer $i$, in its neighborhood $P_i \subset \mathbb{R}^2$ \\
$z_{ij}$ & continuous & quantity transported from $i$ to $j$ (capacity, subroute elimination) \\
$\xi_{ij}$ & continuous & upper bound on the Euclidean distance between $p_i$ and $p_j$ \\
$\gamma_{ij}$ & continuous & distance from $p_i$ to $p_j$ if $i$ precedes $j$, and zero otherwise \\
$w_{i\ell}$ & binary & equals 1 if $p_i$ belongs to the convex polygon $P_{i\ell}$ \\
\hline
\end{tabular}
\caption{Variables of the MIQCP model of the CVRPN.}
\label{tabvar}
\end{table}

The MIQCP model of the CVRPN follows. The model is based on the mixed-integer linear programming (MILP) model presented by~\textcite{Waters1988} for the capacitated vehicle routing problem (CVRP). This formulation serves two purposes. First, it provides a precise mathematical description of the problem addressed in this work. Second, it enables the computation of provably optimal solutions for small-scale instances using a commercial optimization solver. These instances, together with their provably optimal solutions, will later be used as part of the computational experiments to evaluate the effectiveness of the proposed solution method.

{\small\allowdisplaybreaks
\begin{align}
\Minimize & \sum_{i=0}^m \sum_{\substack{j=0\\j\neq i}}^m \gamma_{ij} & \label{cvrpn1} \\
\mbox{s.t. } 
& x_{ij} \in \{0,1\} & i=0,\dots,m, j=0,\dots,m, i\neq j \label{cvrpn11} \\
& \sum_{j = 1}^{m} x_{0j} = \sum_{j = 1}^{m} x_{j0} & \label{cvrpn2} \\
& \sum_{\substack{j=0\\j\neq i}}^{m} x_{ji} = 1 & i = 1,\dots,m\label{cvrpn5} \\
& \sum_{\substack{j=0\\j\neq i}}^{m} x_{ij} = 1 & i = 1,\dots,m\label{cvrpn6} \\
& z_{ij} \leq \mathcal{C} x_{ij} & i=0,\dots,m, j=0,\dots,m, i\neq j \label{cvrpn8} \\ 
& \sum_{\substack{\ell = 0\\ \ell \neq i}}^{m} z_{\ell i} - \sum_{\substack{\ell = 0\\ \ell \neq i}}^{m} z_{i\ell} = d_i & i=1,\dots,m \label{cvrpn7} \\
& \sum_{j = 1}^{m} z_{0j} = \sum_{i=1}^{m} d_i & \label{cvrpn9} \\
& \sum_{i = 1}^{m} z_{i0} = 0 & \label{cvrpn10} \\
& z_{ij} \geq 0 & i=0,\dots,m, j=0,\dots,m, i\neq j \label{cvrpn12} \\
& \xi_{ij}^2 \geq (p^x_i - p^x_j)^2 + ( p^y_i - p^y_j)^2 & i=0,\dots,m, j=0,\dots,m, i\neq j \label{cvrpn13} \\
& \xi_{ij} \geq 0 & i=0,\dots,m, j=0,\dots,m, i\neq j \label{cvrpn14} \\
& \gamma_{ij} \leq M x_{ij} & i=0,\dots,m, j=0,\dots,m, i\neq j \label{cvrpn21} \\
& \gamma_{ij} \geq \xi_{ij} - M (1-x_{ij}) & i=0,\dots,m, j=0,\dots,m, i\neq j \label{cvrpn22} \\
& \gamma_{ij} \geq 0 & i=0,\dots,m, j=0,\dots,m, i\neq j \label{cvrpn23} \\
& a_{i \ell k} p_i^x + b_{i \ell k} p_i^y + c_{i \ell k} \leq 0 & i = 1,\dots,m \mid \nu_i=1, \ell=1, k = 1,\dots,n_{i1} \label{cvrpn24} \\
& a_{i \ell k} p_i^x + b_{i \ell k} p_i^y + c_{i \ell k} \leq M_{i \ell k} (1-w_{i \ell}) & i = 1,\dots,m \mid \nu_i>1, \ell = 1,\dots,\nu_i, k = 1,\dots,n_{i\ell} \label{cvrpn25} \\
& w_{i \ell} \in \{0,1\} & i = 1,\dots,m \mid \nu_i>1, \ell = 1,\dots,\nu_i \label{cvrpn26} \\
& \sum_{\ell=1}^{\nu_i} w_{i \ell} = 1 & i = 1,\dots,m \mid \nu_i>1 \label{cvrpn27}
\end{align}}

The binary variables $x_{ij}$, whose domain is defined by the constraints~\eqref{cvrpn11}, indicate whether customer $i$ immediately precedes customer $j$ on a route ($x_{ij}=1$) or not ($x_{ij}=0$). Using these variables, constraint~\eqref{cvrpn2} states that the same number of vehicles that leave the depot must return to the depot, while constraints~(\ref{cvrpn5},\ref{cvrpn6}) state that each customer is visited exactly once. The continuous variables~$z_{ij}$ correspond to the quantity transported from customer~$i$ to customer~$j$. Constraints~\eqref{cvrpn8} state that this quantity cannot exceed the maximum capacity of the truck and must be zero if customer~$i$ does not precede customer~$j$ on a route. Constraints~\eqref{cvrpn7} state that the quantity transported to customer~$i$ from its predecessor minus the quantity transported from~$i$ to its successor must match the demand of customer~$i$. Constraint~\eqref{cvrpn9} states that the sum of the quantities leaving the depot, in the direction of the first customer on each route, must match the total demand. Constraint~\eqref{cvrpn10} states that nothing should be transported from the last customer on each route back to the depot. Constraints~\eqref{cvrpn12} state that transported quantities must be non-negative. These constraints related to variables~$z_{ij}$ are responsible for the elimination of subroutes. So far, the constraints described are part of the model for the CVRP presented by~\textcite{Waters1988}. Next, we describe the constraints that state that each visit point~$p_i$ must belong to the neighborhood~$P_i$ of customer~$i$ and the constraints that model the distances between consecutive points on a route. 

Constraints~(\ref{cvrpn13},\ref{cvrpn14}) state that $\xi_{ij}$ is an upper bound on the Euclidean distance between points~$p_i$ and~$p_j$. For each pair of $i$ and $j$, constraints~(\ref{cvrpn13},\ref{cvrpn14}) are equivalent to $\xi_{ij} \geq \| p_i - p_j \|_2$. This is a second-order cone (SOC) constraint, which defines a convex feasible set that can be handled by commercial exact solvers such as CPLEX or Gurobi.
An equality in~\eqref{cvrpn13} is avoided because it would define a non-convex feasible set. The inequality is nevertheless sufficient: because the objective function minimizes the sum of the variables $\gamma_{ij}$, at optimal solutions $\xi_{ij}$ coincides with $\| p_i - p_j \|_2$ whenever customer~$i$ precedes customer~$j$ on a route. CPLEX will be used in numerical experiments to find optimal solutions for small instances as a way of evaluating the effectiveness of the proposed method.
These are the constraints that make the model an MIQCP model instead of a plain MILP model. Note that the problem can be formulated as an MILP if the Manhattan or Chebyshev distances are considered instead of the Euclidean distance. Constraints~(\ref{cvrpn21},\ref{cvrpn22},\ref{cvrpn23}) along with the objective function~\eqref{cvrpn1} ensure that, at optimal solutions, $\gamma_{ij}$ coincides with~$\xi_{ij}$ if customer~$i$ precedes customer~$j$ in a route and is zero otherwise. Constraints~\eqref{cvrpn24} state that $p_i \in P_{i}$ and apply only to cases where~$P_i$ is convex, i.e., all~$i$ such that $\nu_i=1$ and~$P_i$ is merely formed by the convex polygon~$P_{i1}$. Constraints~\eqref{cvrpn25} state that $p_i \in P_{i\ell}$ if $w_{i\ell}=1$ and are void if $w_{i\ell}=0$. Constraints~(\ref{cvrpn26},\ref{cvrpn27}) state that the variables $w_{i\ell}$ are binary and that, for each fixed~$i$, they must add up to~1. This implies that, when~$P_i$ is formed by the union of~$\nu_i$ disjoint convex polygons~$P_{i\ell}$, $p_i$ must belong to exactly one of them. The objective function~\eqref{cvrpn1} minimizes the sum of the distances between the visit points of consecutive customers on the routes, including the distances from the depot to the first and last visit points on each route.

In constraints~(\ref{cvrpn21},\ref{cvrpn22}), $M>0$ is a sufficiently large constant. Consider that all polygons~$P_i$ and the depot~$p_0$ are contained in a rectangle $[\underline{p}^x,\overline{p}^x] \times [\underline{p}^y,\overline{p}^y]$. In~(\ref{cvrpn21},\ref{cvrpn22}), it would be enough to take $M$ as $\| (\overline{p}^x,\overline{p}^y) - (\underline{p}^x,\underline{p}^y) \|$, which is an upper bound on the distance between any two visit points. In~\eqref{cvrpn25}, the constant $M_{i \ell k} \geq 0$ must be such that the constraint is void if $w_{i\ell}=0$, i.e., it must be an upper bound on the value that $a_{i \ell k} p_i^x + b_{i \ell k} p_i^y + c_{i \ell k}$ may assume. Because $p_i \in P_i$ in every feasible solution, it suffices to take
\[
M_{i \ell k} = \max_{(p^x,p^y)^T \in P_i} \left\{ a_{i \ell k} p^x + b_{i \ell k} p^y + c_{i \ell k} \right\}
= \max_{(v^x,v^y)^T \in \mathcal{V}(P_i)} \left\{ a_{i \ell k} v^x + b_{i \ell k} v^y + c_{i \ell k} \right\},
\]
where $\mathcal{V}(P_i)$ denotes the set of vertices of~$P_i$. The second equality holds because a linear function attains its maximum over a polygon at one of its vertices; so each constant is computed by simply evaluating a linear function at the vertices of~$P_i$. Note that this computation requires the vertices of each polygon~$P_i$ to be known, which is the case here, because the neighborhoods are polygons given by their vertices. These constants are never larger than, and are usually much smaller than, a single constant valid for all constraints~\eqref{cvrpn25}. They do not change the feasible set of the model, but they yield a tighter continuous relaxation. Regarding numerical stability, in our implementation the coordinates are translated to place the depot at the origin and, if necessary, divided by a common factor to fit in $[-1,1]^2$, and the coefficients of each constraint~(\ref{cvrpn24},\ref{cvrpn25}) are normalized so that $|a_{i \ell k}| + |b_{i \ell k}| + |c_{i \ell k}| = 1$. Consequently, the constants $M$ and $M_{i \ell k}$ have the same order of magnitude as the coordinates, which avoids the numerical difficulties associated with large big-$M$ values. All reported values refer to the original coordinates.

\section{Computing visit points with block coordinate descent} \label{sec3}

Given a visit sequence, i.e., a sequence that starts at the depot, visits customers $i_1, i_2, \dots, i_{\kappa}$, and ends at the depot, the problem of finding the visit point for each of the customers on the sequence can be written as problem~\eqref{proba}. The MIQCP model for this problem is analogous to model \eqref{cvrpn1}--\eqref{cvrpn27}, but with variables $x_{ij}$ and $z_{ij}$ fixed and restricted only to the route at hand. As with \eqref{cvrpn1}--\eqref{cvrpn27}, the optimal solution to this problem may be obtained using a commercial solver such as CPLEX or Gurobi. In this paper, we propose solving it using a block coordinate descent (BCD) method. What follows is a description of a version of BCD for problem (1). This version applies a cyclic selection of blocks, with each block corresponding to a visit point. For a more general version, which includes other ways of identifying blocks and selecting the block to be used in each iteration, see, for example, \textcite{Birgin2022}. A similar scheme, which moves one visit point at a time, was used by~\textcite{Gao2021}, without a description of how each visit point is computed, and by~\textcite{Kulich2022}, who minimize over the edges of each polygon using a rubber-band algorithm originally proposed for pairwise disjoint polygons. Unlike the latter, the subproblem solved here also considers the interior of the polygon, where the optimal point may lie when the segment joining the two neighboring visit points intersects the polygon (Section~\ref{sec31}).

To simplify the presentation, we denote $\pi_j \equiv p_{i_j}$ and $\mathcal{P}_j \equiv P_{i_j}$ for $j=1,\dots,\kappa$ and restate problem~\eqref{proba} as
\begin{equation} \label{probb}
\Minimize_{\{ \pi_j \in \mathcal{P}_j, j=1,\dots,\kappa \}} \Psi(\pi_1,\dots,\pi_{\kappa}) := \sum_{j=0}^\kappa \| \pi_{j+1} - \pi_j \|,
\end{equation}
where, as in~\eqref{proba}, $\pi_0$ and $\pi_{\kappa+1}$ refer to the known geographical position of the depot. A complete cycle of the BCD method consists of varying the index $\ell$ from 1 to $\kappa$ and, for each $\ell$, solving
\begin{equation} \label{probc}
\pi_{\ell} = \argmin_{\{ \pi \in \mathcal{P}_{\ell} \}} \left\{ \| \pi_{\ell+1} - \pi \| + \| \pi - \pi_{\ell-1} \| \right\}.
\end{equation}
The subproblem in~\eqref{probc} is derived from~\eqref{probb} by considering only $\pi_{\ell}$ as a variable while fixing all other $\pi_j$ with $j \neq \ell$, and omitting the constant terms---those not depending on $\pi_{\ell}$---in the objective function. The minimization variable is denoted by $\pi$, and $\pi_{\ell}$ is redefined as an optimal solution of the subproblem. The type of subproblem that appears in~\eqref{probc} is very simple, contains only two variables, and its optimal solution can be efficiently computed. At the end of each cycle, we check whether a stopping criterion is satisfied. If the stopping criterion is satisfied, the visit points are given by $\pi_j$ for $j=1,\dots,\kappa$. If the stopping criterion is not satisfied, a new cycle is started. In practice, we consider two stopping criteria related to the maximum number of cycles of the BCD method and lack of progress in the value of $\Psi$ in consecutive iterations, namely, 1000 cycles and a relative reduction of $\Psi$ smaller than $1.5\times 10^{-8}$ in a cycle. We show in Section~\ref{sec31} how the subproblem in~\eqref{probc} is solved for a fixed $\ell$. When the minimizer of a subproblem is not unique, the current visit point is kept unless another candidate is strictly better, and ties among candidates are broken by the order of the vertices of the polygon.

Before proceeding, it is important to mention that the algorithm described above is a simplification of Algorithm~3.1 presented by~\textcite{Birgin2022}. In Algorithm~3.1 of~\textcite{Birgin2022}, after solving \eqref{probc}, a decrease condition is tested. If this condition is not satisfied, one or more subproblems similar to~\eqref{probc}, augmented by the regularization term $\sigma \| \pi - \pi_{\ell} \|^2$ for increasing values of the regularization parameter $\sigma>0$, must be solved. These extra steps are essential for obtaining worst-case complexity results of the method and ensuring convergence to stationary points, as presented by~\textcite{Birgin2022}. Without these extra steps, there is a possibility that the method will cycle; see the discussion in~\textcite{Amaral2022}. However, solving subproblems of type~\eqref{probc} augmented with the regularization term is not as simple as solving~\eqref{probc}. While the solution of~\eqref{probc} can be performed in a finite number of steps, the solution of its regularized variant would require a more expensive iterative method with asymptotic convergence. Considering that the cases in which the BCD method cycles are pathological and extremely rare in practice, and that we are using BCD within the context of a heuristic method, we believe that the basic version used here is the most appropriate option. The main purpose of the proposed BCD procedure is not to guarantee global optimality of the visit-point problem, but rather to provide high-quality route evaluations at a computational cost compatible with large-scale neighborhood search. Accordingly, the solution computed by BCD may depend on the initial visit points and carries no guarantee of global optimality: in this paper, CPLEX plays the role of an exact method, whereas BCD, for the visit points of a fixed sequence, and GVNS combined with BCD, for the complete CVRPN, are heuristics.

\subsection{Solving the BCD subproblems} \label{sec31}

\begin{figure}[ht!]
\centering
\includegraphics[scale=1.0]{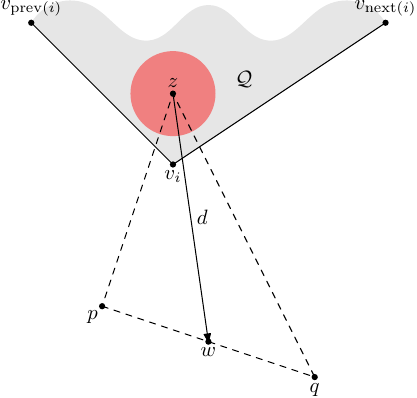}
\caption{This figure illustrates that a point $z$ in the interior of $\mathcal{Q}$ cannot be a solution of~\eqref{probd} when the segment $[p,q]$ does not intersect the polygon $\mathcal{Q}$. This follows from the fact that for any $w$ in the segment $[p,q]$, the direction $d=w-z$ is both a feasible and a descent direction from $z$.}
\label{fig2}
\end{figure}

In this section, we show how to solve problem~\eqref{probc} for a fixed $\ell \in \{1,\dots,\kappa\}$. To simplify the presentation, we restate the problem as 
\begin{equation} \label{probd}
\Minimize_{z \in \mathcal{Q}} \phi(z) := \| p - z \| + \| z - q \|,
\end{equation}
where $p, q \in \mathbb{R}^2$ and $\mathcal{Q}$ is a polygon, possibly non-convex. The procedure assumes that $\mathcal{Q}$ is a simple polygon whose consecutive vertices are distinct, so that no edge is degenerate. By the triangle inequality, the optimal value of~$\phi$ is greater than or equal to $\|p-q\|$. If $[p,q] \cap \mathcal{Q} \neq \emptyset$, then any point $z^{\star} \in [p,q] \cap \mathcal{Q}$ achieves the optimal value and is therefore a solution. This includes the case in which $p$ or $q$ lies on the boundary of $\mathcal{Q}$. We now intend to show that, when $[p,q] \cap \mathcal{Q} = \emptyset$, a solution $z^{\star}$ is necessarily on the boundary of $\mathcal{Q}$. Let $z$ be in the interior of $\mathcal{Q}$. Because $[p,q] \cap \mathcal{Q} = \emptyset$, we have $z \notin [p,q]$ and, hence, $\phi(z) > \|p-q\|$, because equality holds in the triangle inequality if and only if $z \in [p,q]$. This holds even if $p$, $q$, and $z$ are collinear or if $p=q$. Let $w$ be a point in the segment $[p,q]$ and let $d=w-z$; see Figure~\ref{fig2}. Because $z$ is in the interior of $\mathcal{Q}$, there exists $\bar \epsilon \in (0,1]$ such that $z+\epsilon d \in \mathcal{Q}$ for all $\epsilon \in (0,\bar \epsilon]$. Moreover, because $\phi$ is convex and $\phi(w)=\|p-q\|<\phi(z)$, we have $\phi(z+\epsilon d) \leq (1-\epsilon)\phi(z)+\epsilon\phi(w) < \phi(z)$ for all $\epsilon \in (0,\bar \epsilon]$. Therefore, $z$ does not satisfy the definition of local minimizer and cannot, thus, be a solution to~\eqref{probd}. We conclude that, to solve~\eqref{probd} when $[p,q] \cap \mathcal{Q} = \emptyset$, we need only to search for a solution on the boundary of~$\mathcal{Q}$. The above argument does not depend at all on $\mathcal{Q}$ being a polygon and remains valid for a customer neighborhood given by any compact set. Next, we describe how to find a solution on the boundary of $\mathcal{Q}$. This description depends heavily on the fact that $\mathcal{Q}$ is a polygon.

\begin{figure}[ht!]
\centering
\resizebox{\textwidth}{!}{
\begin{tabular}{ccc}
\includegraphics[scale=1.0]{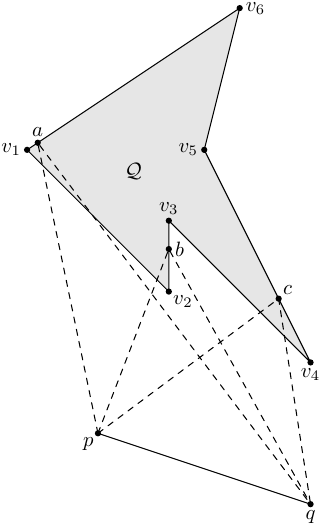} & &
\includegraphics[scale=1.0]{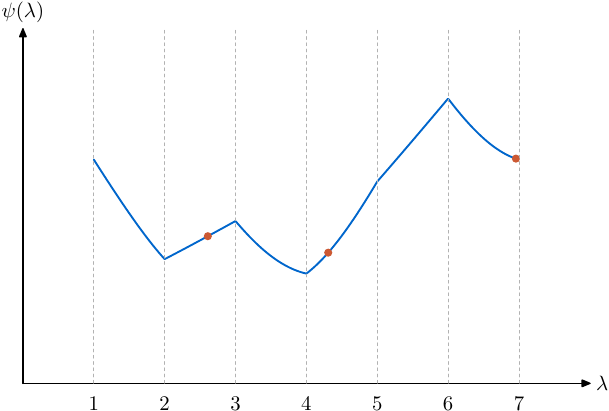}
\end{tabular}}
\caption{(a) Graphical representation of the parameterization of the boundary of a non-convex polygon $\mathcal{Q}$ with vertices $v_1, v_2, \dots, v_6$ in counterclockwise order. The whole boundary corresponds to varying $\lambda$ between 1 and 7. The point $a=v_6 + 0.95 (v_1-v_6)$ corresponds to $\lambda=6.95$; the point $b=v_2+0.6(v_3-v_2)$ corresponds to $\lambda=2.6$; the point $c=v_4+0.3(v_5-v_4)$ corresponds to $\lambda=4.3$. In this example, $[p,q] \cap \mathcal{Q}=\emptyset$. (b) Graph of $\psi(\lambda)$. The graph highlights in red the pairs $(\lambda,\psi(\lambda))$ that correspond to points $a$, $b$, and $c$. It is clear that the function is convex and differentiable between integer values of $\lambda$ and non-differentiable for integer values of $\lambda$. A simple inspection of both figures shows that $z^{\star}=z(4)=v_4$.}
\label{fig3}
\end{figure}

We now outline how to compute a solution to~\eqref{probd} on the boundary of $\mathcal{Q}$ when $[p,q] \cap \mathcal{Q} = \emptyset$. Recall that we are assuming that $\mathcal{Q}$ has $h \geq 3$ vertices\footnote{Solving~\eqref{probd} when $\mathcal{Q}$ has a single vertex is trivial, while the case with two vertices is similar to the case with three or more vertices, but only two vertices and one edge need to be verified. Treating the case with two vertices in the same way as situations involving three or more vertices would result in the same edge being tested twice.} that here we name $v_{\ell}$ ($\ell=1,\dots,h$) in counterclockwise order. To simplify the presentation, given $\ell \in \{1,\dots,h\}$, we define $\mathrm{next}(\ell)=\ell+1$ whenever $\ell<h$ and $\mathrm{next}(h)=1$. Analogously, we define $\mathrm{prev}(\ell)=\ell-1$ whenever $\ell>1$ and $\mathrm{prev}(1)=h$. Because we limit the search to the boundary of $\mathcal{Q}$, we can parameterize its boundary with a single scalar parameter~$\lambda$. Intuitively, $\lambda$ describes a walk along the boundary of $\mathcal{Q}$ in counterclockwise order, starting at vertex~$v_1$: the integer part of $\lambda$ identifies the edge being traversed, while its fractional part indicates how far along that edge the point $z(\lambda)$ lies. This allows us to rewrite~\eqref{probd} as
\begin{equation} \label{probe}
\Minimize_{\lambda \in \mathbb{R}} \psi(\lambda) := \| p - z(\lambda) \| + \| z(\lambda) - q \|
\mbox{ subject to } 1 \leq \lambda \leq h+1,
\end{equation}
where, for $i=1,\dots,h$, if $\lambda \in [i,i+1]$, then $z(\lambda) := v_i + ( \lambda - i) ( v_{\mathrm{next}(i)} - v_i )$. Note that if $\lambda=i$ then $z(\lambda)=v_i$ for $i=1,\dots,h$ and that non-integer values of $\lambda$ correspond to points in the interior of the edges of $\mathcal{Q}$. Concretely, if $\lambda \in (i,i+1)$, then $z(\lambda) \in (v_i,v_{\mathrm{next}(i)})$. The key to finding a solution to~\eqref{probe} is that the function $\psi(\lambda)$ is convex (sum of norms of affine functions) within each interval $\lambda \in [i,i+1]$. It is strictly convex in this interval unless both $p$ and $q$ lie on the straight line containing the edge $[v_i,v_{\mathrm{next}(i)}]$. In that case, because $[p,q] \cap \mathcal{Q} = \emptyset$, $p$ and $q$ lie on the same side of the edge along that line, $\psi$ is affine and non-constant in $[i,i+1]$, and its minimum in this interval is attained at a vertex. The fact that $[p,q] \cap \mathcal{Q} = \emptyset$ guarantees that $\psi$ is differentiable in $(i,i+1)$. See Figure~\ref{fig3}. Therefore, $z(\lambda) \in [v_i,v_{\mathrm{next}(i)}]$ can be a solution only if $\lambda \in [i,i+1]$ is such that $\lambda=i$, $\lambda=i+1$ or $\lambda \in (i,i+1)$ and $\psi'(\lambda)=0$. In other words, the candidate solutions are the vertices of $\mathcal{Q}$ and the points in the interior of the edges that nullify the derivative of $\psi$. Among these candidates, the one with the lowest value of $\psi$ is the desired solution.

For $i=1,\dots,h$, if $\lambda \in (i,i+1)$, the derivative of $\psi$ is given by
\begin{equation} \label{probeder}
\psi'(\lambda) = \left\langle - \frac{p-z(\lambda)} {\|p-z(\lambda)\|} + \frac{z(\lambda)-q} {\|z(\lambda)-q\|}, v_{\mathrm{next}(i)}-v_i \right\rangle;
\end{equation}
and therefore $\psi'(\lambda)=0$ if and only if
\begin{equation} \label {nullpsider2}
\frac{\langle z(\lambda)-q, v_{\mathrm{next}(i)}-v_i \rangle}{\|z(\lambda)-q\|} = 
\frac{\langle p-z(\lambda), v_{\mathrm{next}(i)}-v_i \rangle}{\|p-z(\lambda)\|}.
\end{equation}
By squaring both sides (to eliminate the square root) and cross-multiplying, we obtain a quartic equation, which reduces to the quadratic equation given by $A t^2 + B t + C = 0$ with $A = \delta_q^2-\delta_p^2$, $B = -2 (\alpha_p \delta_q^2 - \alpha_q \delta_p^2)$, and $C = \alpha_p^2 \delta_q^2 - \alpha_q^2 \delta_p^2$, where $t=\lambda-i$ and, for $r \in \{p,q\}$, $\alpha_r=\langle r-v_i,d_i\rangle/\|d_i\|^2$ and $\delta_r=\|r-v_i-\alpha_r d_i\|$, with $d_i=v_{\mathrm{next}(i)}-v_i$. When $\delta_p=\delta_q>0$, the equation is linear and its solution is $t^{\star}=(\alpha_p+\alpha_q)/2$. Otherwise, because the quadratic equation is equivalent to $\delta_q^2 (t-\alpha_p)^2 = \delta_p^2 (t-\alpha_q)^2$, its roots are given by $t_1=(\alpha_p \delta_q + \alpha_q \delta_p)/(\delta_p+\delta_q)$ and $t_2=(\alpha_p \delta_q - \alpha_q \delta_p)/(\delta_q-\delta_p)$. Both roots are well defined, $t_2$ because $\delta_p \neq \delta_q$, and $t_1$ because $\delta_p+\delta_q>0$, as the case $\delta_p+\delta_q=0$, in which both $p$ and $q$ lie on the line containing the edge $[v_i,v_{\mathrm{next}(i)}]$, was treated in the previous paragraphs. We substitute $t_1$ and $t_2$ into~\eqref{nullpsider2} and define $t^{\star}$ as the one that satisfies it. This root corresponds to the minimizer of $\psi$ on the line passing through the vertices $v_i$ and $v_{\mathrm{next}(i)}$. Let $\lambda_i^{\star}=i+t^{\star}$. If $\lambda_i^{\star} \in (i,i+1)$, then $z(\lambda_i^{\star}) \in (v_i,v_{\mathrm{next}(i)})$ and $z(\lambda_i^{\star})$ is a candidate solution to~\eqref{probd}. Let $\mathcal{S}$ be the set of candidates $z(\lambda_i^{\star})$ plus the vertices of $\mathcal{Q}$. The solution $z^{\star}$ to~\eqref{probd} is given by $z^{\star} = \argmin_{\{ z \in \mathcal{S}\}} \left\{ \phi(z) \right\}$.

\section{General variable neighborhood search} \label{sec4}

General variable neighborhood search (GVNS) is a variation of the variable neighborhood search \parencite[VNS;][]{Hansen2001,Hansen2009,Mladenovi1997} metaheuristic designed to enhance exploration and intensification in the search process. VNS systematically explores different search neighborhoods to escape local optima by alternating between a shaking phase (to diversify) and a local search phase (to improve the solution). However, it typically uses a single neighborhood structure in each local search step. GVNS strengthens this approach by incorporating a variable neighborhood descent (VND) within the local search phase. The VND method applies several neighborhood structures sequentially. This deeper, more exhaustive search yields a better-quality solution before returning to the shaking step. This hierarchical structure makes GVNS more robust and effective than the basic VNS, especially for complex combinatorial problems, such as the VRP and its variants.

When applied to VRP variants, GVNS has delivered strong performance. Recent works~\parencite{Bezerra2018,deArmas2015,Kytjoki2007,Masmoudi2025,Mladenovi2012,Tadaros2024,Todosijevi2017} show that GVNS can achieve near-optimal or optimal solutions for many VRP benchmark instances, performing competitively against leading metaheuristics such as Tabu Search, Genetic Algorithms, and Ant Colony Optimization. The VRP's structure allows, in a natural fashion, for multiple neighborhood definitions, primarily based on route swaps and customer relocations. This variety of neighborhood definitions fits perfectly within the GVNS framework, enabling it to effectively navigate the VRP's complex solution space. In this study, we considered the intra-route neighborhoods \textbf{(i)} 2-Opt, \textbf{(ii)} 3-Opt, and \textbf{(iii)} Or-Opt and the inter-route neighborhoods \textbf{(iv)} 1-0-Exchange, \textbf{(v)} 1-1-Exchange, \textbf{(vi)} 1-2-Exchange, \textbf{(vii)} Relocate, \textbf{(viii)} Cross-Exchange, and \textbf{(ix)} 2-Opt-Star~\parencite[see][]{aarts2018,TothVigo2002VRP}. 

In the following, we describe implementation details of each of the considered neighborhoods. We implemented the 2-Opt neighborhood as described in the literature. The 3-Opt neighborhood, as described in the literature, includes the movements of 2-Opt. In our implementation, as we intend to use both combined, we implemented in 3-Opt only those movements that are not 2-Opt. The Or-Opt neighborhood, given a route, considers segments with $r$ consecutive customers, with $r$ varying from $r_{\min}$ to $r_{\max}$, and reallocates the segment in all possible positions on the route. Unlike the classic version described in the literature, when reallocating, we also consider reallocating the inverted sequence. We implemented the 1-0-Exchange and the 1-1-Exchange neighborhoods exactly as described in the literature. The 1-2-Exchange considers two routes; from one it takes a segment with a single customer, from the other it takes a segment with two consecutive customers, and exchanges these segments on the routes. Unlike the classic version described in the literature, when inserting the segment with two customers, we consider two possibilities: the usual insertion and the insertion with the inverted segment. The Relocate neighborhood is an extension of 1-0-Exchange. It takes segments with $r$ consecutive customers from one route, with $r$ varying from $r_{\min}$ to $r_{\max}$, removes them from the route, and inserts them into another route. Unlike the classic version described in the literature, when inserting, we also consider the inverted sequence insertion. The Cross-Exchange neighborhood is an extension of 1-1-Exchange. It considers two routes and exchanges consecutive segments with $r_1$ and $r_2$ customers from each route. Unlike the classic version described in the literature, we consider the four combinations in which the segments are reinserted, inverted or not. In the experiments, we consider all pairs $(r_1, r_2)$ in a predefined set of pairs $R$. The pairs $(1,0)$ and $(0,1)$ would generate the movements of the 1-0-Exchange neighborhood. Other pairs containing~0, either as the first or second coordinate, would generate movements in the Relocate neighborhood. The pair $(1,1)$ would generate movements in the 1-1-Exchange neighborhood, and the pairs $(1,2)$ and $(2,1)$ would generate movements in the 1-2-Exchange neighborhood. Therefore, the set $R$ can be defined as $R = \{ (r_1,r_2) \;|\; r_{\min} \leq r_1, r_2 \leq r_{\max}\} \setminus \{ (1,1), (1,2), (2,1) \}$, with $r_{\max} \geq r_{\min} \geq 1$. We implemented 2-Opt-Star exactly as described in the literature. 

Because there is no fixed cost associated with the use of a vehicle, the number of routes is not fixed in advance. To allow GVNS to explore solutions with different numbers of routes, a solution is always represented with $m$ routes, some of which may be empty. In the 1-0-Exchange and Relocate neighborhoods, the destination route may be an empty route, so that these neighborhoods can open a new route with a single customer (1-0-Exchange) or with a segment of $r$ consecutive customers, $r_{\min} \leq r \leq r_{\max}$ (Relocate). Because all empty routes are equivalent, only one of them is considered as a destination route in each exploration of these neighborhoods. Conversely, a route becomes empty, and is thus closed, whenever all its customers are moved to other routes. The remaining inter-route neighborhoods (1-1-Exchange, 1-2-Exchange, Cross-Exchange, and 2-Opt-Star) operate on pairs of non-empty routes. In the shaking step, empty routes may also be selected as destination routes.

The natural value for the parameter $r_{\min}$ in the Or-Opt and Cross-Exchange neighborhoods is $r_{\min}=1$. In the case of the Relocate neighborhood, the natural value is $r_{\min}=2$, because with $r=1$ the movements are those already generated by the 1-0-Exchange neighborhood. There are no strong reasons for the three neighborhoods to use different values for $r_{\max} \geq 2$. Therefore, we will consider a unified value for all three. This value will be the subject of numerical experimentation. During the parameter calibration phase, we will also evaluate the inclusion of movements with inverted sequences. Inverting a segment may yield a better move, because the best orientation of a segment depends on its new neighbors, but it doubles the number of movements to be examined. These movements will be prohibited if \(\textsc{flip}=false\) and allowed if \(\textsc{flip}=true\). The decision on whether to use one or the other will be the subject of numerical experimentation. Related to the use of neighborhoods, as is usual in local searches, we will analyze whether we should use the first-improvement strategy and abandon the construction of the neighborhood when a neighbor better than the current solution is found, or whether the inspection of the neighborhood should be complete in order to choose the best neighbor as the new solution. We will consider the parameter \textsc{otype} during the calibration phase, assigning value F to the first-improvement strategy and value B to the best-improvement strategy. When the first-improvement strategy is used, the routes can be examined either in a fixed order or in a new random order each time the neighborhood is explored, so that no route is systematically examined first. This option, which we will refer to as $\textsc{randroutes} \in \{true,false\}$, will also be subject to calibration.

Our implementation of GVNS follows Algorithms~4, 5, and~8 presented on pages~372, 373, and~376 of~\textcite{Hansen2009}. It is interesting to note that the method has no parameters other than those related to the neighborhoods. The neighborhoods considered are the nine neighborhoods mentioned above. There is a shake associated with each neighborhood, and the shake consists of making a random movement of the associated neighborhood. However, preliminary experiments have shown that a single movement of a single neighborhood does not provide the necessary diversification to the method. We therefore consider a shake that performs $q_{\max}$ consecutive random movements, each in a randomly selected neighborhood. The quantity $q_{\max}$ of consecutive random movements that make up the shake will be subject to numerical calibration. Another factor relevant to the performance of the method is the order in which the neighborhoods are considered in the VND. In the numerical experiments, we consider four basic possibilities. The first is to consider the neighborhoods in the order in which they were described. The second possibility is to consider first the six inter-route neighborhoods (in the order in which they were described) and then the three intra-route neighborhoods (also in the order in which they were described). In these two options, the order is maintained throughout the entire execution of the method. The third and fourth possibilities consist of starting with the first and the second order, respectively, and changing the order to a random order each time the method completes a cycle of nine shakes without improvement of the incumbent. The choice of this parameter, which we will refer to as $\text{\textsc{neighborhood\_order}} \in \{1,2,3,4\}$, will be subject to calibration. The VND first explores the nine neighborhoods with the visit points of the current solution kept fixed, then with the visit points restricted to the vertices of the polygons and computed by a shortest-path algorithm, and finally with the visit points computed by BCD, starting from those of the current solution, so that BCD is used only when the cheaper evaluations no longer find improving moves.

\section{Numerical experiments} \label{sec5}

In this section, we present numerical experiments to evaluate the performance of the proposed metaheuristic enhanced with block coordinate descent. Different aspects of the proposed model and method will be individually evaluated. Section~\ref{sec51} describes two simple constructive heuristics that deliver a solution to be used as the initial solution for either the GVNS metaheuristic or an exact commercial solver. In Section~\ref{sec52}, we introduce small- and medium-scale CVRPN instances, each with 3 to 60 customers. These instances are used, in Section~\ref{sec54}, to compare two forms of calculating the visit points: using BCD versus solving an MIQCP model with the commercial solver CPLEX. The evaluation focuses on the effectiveness and efficiency aspects of BCD. In Section~\ref{sec55}, we rely on instances from the capacitated VRP literature to construct large-scale instances of the CVRPN with polygonal customer neighborhoods. Section~\ref{sec56} describes the calibration of the parameters of the proposed GVNS metaheuristic. Sections~\ref{sec57} and~\ref{sec58} describe the performance of the GVNS metaheuristic in the small-, medium-, and large-scale instances of the CVRPN. Table~\ref{tabexp} summarizes the experimental design.

\begin{table}[ht!]
\centering
\setlength{\tabcolsep}{5pt}
\begin{tabular}{lll}
\hline
 & Small- and medium-scale & Large-scale \\
\hline
Number of instances & 134 & 318 (106 $\times$ 3 areas) \\
Number of customers & 3 to 60 & 15 to 385 \\
Neighborhoods & 3 to 12 vertices, no overlaps & 3 to 8 vertices, overlaps allowed \\
Origin & randomly generated & adapted from CVRPLib \\
Exact solver & CPLEX, one hour & CPLEX, one hour \\
GVNS & 10 runs of 60 or 1200 seconds & 10 runs of 1200 seconds \\
Long GVNS run & -- & one run of one or five hours \\
Compared with & CPLEX & CPLEX; CVRP-based upper bound \\
\hline
\end{tabular}
\caption{Summary of the experimental design.}
\label{tabexp}
\end{table}

The GVNS metaheuristic was implemented in Fortran~90. To facilitate reproducibility, the complete source code, all generated benchmark instances, and the solutions reported in this paper are publicly available at \url{https://doi.org/10.5281/zenodo.24095379} (see also \url{http://www.ime.usp.br/~egbirgin/}). The experiments were carried out on an Intel i9-12900K (12th Gen) processor operating at 5.2~GHz with 128 GB of RAM running Ubuntu 23.04. The code was compiled by the GNU Fortran compiler GCC (version 12.3.0) with the -O3 optimization directive enabled. Models were solved using IBM ILOG CPLEX Optimization Studio 22.1.2 with its default parameters, except for the number of threads, which was set to one. Because the GVNS metaheuristic (including BCD) is a sequential code, all methods were run on a single thread. All time limits and reported times refer to CPU time. When applied to the complete model, CPLEX starts from the better of the two solutions constructed with the heuristics described in Section~\ref{sec51}.

\subsection{Heuristics to construct feasible solutions} \label{sec51}

To both initialize the GVNS metaheuristic and assist in solving the model through the application of an exact method, it is useful to have a constructive heuristic that can generate an initial feasible solution. In this work, we consider simple adaptations of two constructive heuristics for the VRP, the savings method of \textcite{Clarke1964} and the nearest neighbor heuristic of \textcite{Flood1956}. In both cases, the heuristics are applied to a CVRP instance derived from the corresponding CVRPN instance. This derived instance is built by assigning to each customer the barycenter of its polygonal neighborhood as its geographical position. We check the capacity constraint before applying each heuristic's main step: merging two routes (savings method), or adding a new customer to an existing route (nearest neighbor). The application of either constructive heuristic yields a set of visit sequences. By applying BCD to each sequence, we determine the visit points and obtain a complete set of routes.

\subsection{Generation of small- and medium-scale instances of the CVRPN} \label{sec52}

\begin{figure}[ht!]
\centering
\begin{tabular}{ccc}
\includegraphics{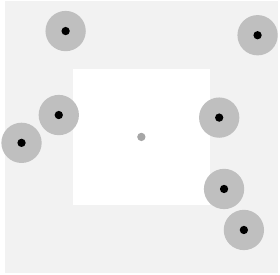} & &
\includegraphics{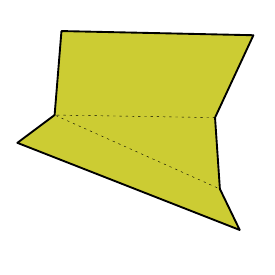} \\
(a) Drawing of the vertices & & (b) Constructed polygon \\
\end{tabular}
\caption{This figure illustrates the construction of a non-convex polygonal customer neighborhood. Figure~(a) shows a random reference point (light gray) and two squares with sides $2 \rho_{\min}$ and $2 \rho_{\max}$ centered on the reference point. It also shows the $m_c=7$ points drawn at random within the larger square and outside the smaller square, with minimum distance $d_{\min}$. Figure~(b) shows the polygon constructed with the drawn points and its division into convex subpolygons. The division is used only in the context of the MIQCP model.}
\label{fig4}
\end{figure}

We generated~134 instances using the scheme proposed by~\textcite{Laporte1984}. The positions of the customers correspond to points that were randomly selected within the square $[-1,1]^2$. However, rather than randomly selecting a single geographical position for each customer~$c$, we first randomly selected a scalar number, say $m_c$. We then randomly selected a cluster of~$m_c$ points that will be associated with customer~$c$. All draws follow a uniform distribution. The generation of points for each cluster considers a random reference point and two squares centered on the reference point with sides of length $2 \rho_{\min}$ and $2 \rho_{\max}$, respectively, where $\rho_{\min}$ and $\rho_{\max}$ are constants that satisfy $0<\rho_{\min}<\rho_{\max} \ll 1$. The cluster points are drawn in such a way that they are outside the smaller square and inside the larger square. Also, points within the same cluster are at a distance $d_{\min} > 0$ or greater from each other, measured in the Euclidean norm. To achieve this, points are drawn within the larger square, discarding those that fall within the smaller square or at a distance of less than $d_{\min}$ from a previously drawn point in the same cluster. See Figure~\ref{fig4}. The reference point is drawn in $[-1+\rho_{\max},1-\rho_{\max}]^2$, so that the entire cluster is contained within $[-1,1]^2$. If the drawn reference point is at a distance (in the infinity norm) of less than $2 \rho_{\max}$ from a reference point of a previous customer, it is discarded, and a new reference point is drawn. The result of this process is a set of ``chunky clusters without overlap''. The depot is always located at the origin.

Up to this point, the problem description could be considered a generalized vehicle routing problem (GVRP) instance. By representing in polar coordinates the points associated with a customer, and ordering them from smallest to largest by angle, we obtain vertices (ordered in a counterclockwise direction) of a non-convex polygon representing the customer's region. In order for the instance to be described with model~\eqref{cvrpn1}--\eqref{cvrpn27}, each non-convex polygon needs to be partitioned into a finite number of convex polygons. We did this using the subroutine \texttt{optimal\_convex\_partition\_2} from the CGAL library~\parencite{cgal:h-pp2-25b}, which computes a partition with the minimum number of convex polygons. For each customer, we randomly selected a demand $d_i$ from the set $\{1, 2, \dots, 10\}$. We generated instances with the number of customers $m$ varying from 3 to 60.

The maximum vehicle capacity~$\mathcal{C}$, the range of the number of points per cluster, and the values of $\rho_{\min}$ and $\rho_{\max}$, all of which depend on the number of customers~$m$, are given in Table~\ref{tabgen}. In all instances, we set $d_{\min} = 0.4 \sqrt{16(\rho_{\max}^2-\rho_{\min}^2)/(\pi m_c)}$. The total number of vertices of the polygons that constitute an instance varies between 9 and 473. The total number of convex polygons that compose the partitions of the non-convex polygonal regions varies from 3 to 140. When applied to the instances, model~\eqref{cvrpn1}--\eqref{cvrpn27} has up to 3,784 binary variables, 11,100 continuous variables, 11,840 linear constraints, and 3,660 quadratic constraints. Tables~1 and~2 in the Supplementary Material show details of the generated instances.

\begin{table}[ht!]
\centering
\begin{tabular}{ccccc}
\hline
Number of customers & $\mathcal{C}$ & Points per cluster & $\rho_{\min}$ & $\rho_{\max}$ \\
\hline
$3 \leq m \leq 15$ & 30 & 3 -- 7 & $0.1\sqrt{3/m}$ & $0.2\sqrt{3/m}$ \\
$16 \leq m \leq 60$ & 50 & 3 -- 12 & $0.05\sqrt{16/m}$ & $0.1\sqrt{16/m}$ \\
\hline
\end{tabular}
\caption{Parameters used to generate the small- and medium-scale instances, for each range of the number of customers~$m$. The columns give the maximum vehicle capacity~$\mathcal{C}$, the range of the number of points per cluster, and the half-sides $\rho_{\min}$ and $\rho_{\max}$ of the squares from which the cluster points are drawn.}
\label{tabgen}
\end{table}

\subsection{Evaluation of using BCD to compute the visit points of a visit sequence} \label{sec54}

Given a visit sequence, i.e., a sequence that starts at the depot, visits customers $i_1, \dots, i_{\kappa}$, and ends at the depot, we propose calculating the visit points using BCD, as described in Section~\ref{sec3}. Alternatively, the visit points can be calculated by solving a simplified version of the MIQCP model~\eqref{cvrpn1}--\eqref{cvrpn27}, which corresponds to fixing the binary variables that define the visit sequence. In this simplified model, only the arcs of the visit sequences, the visit points, the distances between consecutive visit points, and the constraints that require each visit point to belong to its neighborhood are kept, so that the only remaining binary variables are those associated with the convex parts of non-convex neighborhoods. In this section, we aim to show that using BCD is several orders of magnitude more efficient, to the point that without its use, the application of the GVNS metaheuristic proposed in this work would be impractical. For the experiment, we considered the 134 instances introduced in Section~\ref{sec52}. For each of them, we considered the best set of visit sequences (set with the lowest sum of lengths) among the solutions calculated with the constructive heuristics described
in Section~\ref{sec51}. For these visit sequences, we computed the visit points in two ways: by solving the simplified version of the MIQCP model~\eqref{cvrpn1}--\eqref{cvrpn27} and by using BCD. Tables 3 and 4 in the Supplementary Material show the details.

The analysis of the results shows that, in terms of efficiency, BCD is two to four orders of magnitude faster than CPLEX; see Table~\ref{tabbcd}. BCD requires, on average, about three cycles per route (3.1 in the instances with up to 15 customers and 3.2 in the remaining ones), at most 11, and it always stopped because of lack of progress. While BCD uses a CPU time between $10^{-7}$ and $10^{-4}$ seconds, CPLEX takes between $10^{-4}$ and $0.2$ seconds, depending on the instance, and it proved optimality in all 134 instances. Because the CPU time of a single run of BCD is below the resolution of the clock, BCD was run repeatedly until the accumulated CPU time reached at least one second, and its CPU time is this total divided by the number of runs. Of course, the comparison of times should be taken with caution, because we are comparing an exact method with a heuristic method. In terms of the quality of the solution found, CPLEX and BCD find solutions of equivalent quality (relative difference of at most $10^{-5}$) in 120 of the 134 instances. In the remaining 14 instances, BCD finds a slightly better solution, with a relative difference of at most $8.2 \times 10^{-5}$, which is within the optimality tolerance of CPLEX. In summary, in this experiment, BCD reached the same solution quality as CPLEX at a small fraction of its computational cost, which makes the repeated evaluation of routes feasible. Starting from 10 random initial points instead of the barycenters, BCD reached the same route length in 133 of the 134 instances considered in this section.

\begin{table}[ht!]
\centering
\resizebox{\textwidth}{!}{
\begin{tabular}{lr@{\hspace{8pt}}c@{\hspace{8pt}}ccc@{\hspace{8pt}}c@{\hspace{8pt}}ccc@{\hspace{8pt}}c@{\hspace{8pt}}c@{\hspace{8pt}}c@{\hspace{8pt}}c}
\hline
 & & & \multicolumn{3}{c}{CPLEX CPU time (s)} & & \multicolumn{3}{c}{BCD CPU time (s)} & & Ratio CPLEX/BCD & & BCD \#cycles \\
\cline{4-6}\cline{8-10}\cline{12-12}\cline{14-14}
Instances & \# & & min & median & max & & min & median & max & & min -- max & & average \\
\hline
$m \leq 15$ & 65 & & $3.1\times 10^{-4}$ & $2.0\times 10^{-3}$ & $4.3\times 10^{-3}$ & & $4.8\times 10^{-7}$ & $1.8\times 10^{-6}$ & $5.0\times 10^{-6}$ & & $2.0\times 10^{2}$ -- $2.3\times 10^{3}$ & & 6.6 \\
$m > 15$ & 69 & & $4.8\times 10^{-3}$ & $1.6\times 10^{-2}$ & $1.9\times 10^{-1}$ & & $3.4\times 10^{-6}$ & $1.2\times 10^{-5}$ & $2.7\times 10^{-5}$ & & $6.4\times 10^{2}$ -- $1.5\times 10^{4}$ & & 15.4 \\
\hline
All & 134 & & $3.1\times 10^{-4}$ & $5.3\times 10^{-3}$ & $1.9\times 10^{-1}$ & & $4.8\times 10^{-7}$ & $5.1\times 10^{-6}$ & $2.7\times 10^{-5}$ & & $2.0\times 10^{2}$ -- $1.5\times 10^{4}$ & & 11.2 \\
\hline
\end{tabular}}
\caption{Summary of the CPU times required to compute the visit points of a given set of visit sequences using CPLEX and BCD, ratio between the CPU times of CPLEX and BCD (computed for each instance), and average number of BCD cycles per instance (added over all its routes), for the 134 small- and medium-scale instances ($m$ denotes the number of customers). Detailed results are reported in Tables~3 and~4 in the Supplementary Material.}
\label{tabbcd}
\end{table}

\subsection{Large-scale instances derived from literature CVRP instances} \label{sec55}

We consider the CVRP instances available in the CVRPLib library~\parencite{CVRPLib}. The steps for converting a CVRP instance to a CVRPN instance are as follows. The neighborhood of each customer is generated with a procedure similar to the one described in Section~\ref{sec52}, with two differences: the two squares are centered at the customer's geographical coordinates (instead of at a random reference point), and no condition is imposed to prevent the neighborhoods of nearby customers from overlapping. More precisely, for each customer, we randomly select a number of points $nv$ between $nv_{\min}=3$ and $nv_{\max}=8$ and draw $nv$ points inside the square with sides of length $2\rho_{\max}$ and outside the square with sides of length $2\rho_{\min}$, both centered at the customer, discarding points at a distance of less than $0.5\sqrt{16(\rho_{\max}^2-\rho_{\min}^2)/(\pi\, nv)}$ from a previously drawn point of the same customer. The values of $\rho_{\min}$ and $\rho_{\max}=2\rho_{\min}$ are chosen so that the larger square of each customer has an area equal to 0.1\%, 0.5\%, or 1.0\% of the area of the smallest rectangle containing the depot and the customers. We express the cluster points in polar coordinates and sort them in ascending order by angle. The points then represent the vertices of a potentially non-convex polygonal customer neighborhood. We then check if the customer belongs to its neighborhood. If not, we replace a random point in the cluster with the customer, reorder the cluster points by angle, and obtain a new neighborhood that is guaranteed to include the customer. The conversion process preserves the location of the depot. As a result, any solution to the original CVRP instance (in particular, the optimal or best solution reported in the literature) is a feasible solution to the generated CVRPN instance. Consequently, the best-known solution for the CVRP instance can serve as an upper bound for the optimal solution of the CVRPN instance, which we aim to either achieve or improve upon in our numerical experiments.

Using the procedure described above, we generated CVRPN instances based on the A, B, E, F, M, P, and Tai benchmark sets. From each CVRP instance, we generated three CVRPN instances, one for each of the three areas, for a total of 318 CVRPN instances. The three instances derived from the same CVRP instance have the same number of vertices in each neighborhood and differ only in their coordinates. Across these instances, the number of customers ranges from 15 to 385, while the total number of customer neighborhood vertices ranges from 86 to 2,088, and 47\% of the customers have a non-convex neighborhood. Tables~5 and~6 in the Supplementary Material show details of the generated instances. For the largest generated instance, model~\eqref{cvrpn1}--\eqref{cvrpn27} comprises approximately 149,000 binary variables, 447,000 continuous variables, 450,000 linear constraints, and 149,000 quadratic constraints. Because no condition is imposed to prevent overlaps, 70, 95, and 100 of the 106 instances with areas of 0.1\%, 0.5\%, and 1.0\%, respectively, contain at least one pair of customers with overlapping neighborhoods. With an area of 0.1\%, overlaps are particularly frequent in sets B, F, M, and Tai, in which between 11\% and 85\% of the customers of each instance have a neighborhood that overlaps with the neighborhood of at least one other customer. With an area of 1.0\%, this fraction ranges from 80\% to 100\% in these sets. These instances are handled by the proposed method without any modification. 

\subsection{Calibration of the parameters of the GVNS metaheuristic} \label{sec56}

We automatically calibrated the parameters of the GVNS heuristic using Irace~\parencite{irace}. For this purpose, we generated 15 calibration instances with 80 to 200 customers, following the scheme of~\textcite{Uchoa2017}, and converted them into CVRPN instances as described in Section~\ref{sec55}, with neighborhood areas of 0.1\%, 0.5\%, and 1.0\% (five instances each). These instances are not used in the remaining experiments. The search space included the parameters $\text{\textsc{otype}} \in \{\mathrm{B},\mathrm{F}\}$, $r_{\max} \in \{2,\dots,10\}$, $\text{\textsc{neighborhood\_order}} \in \{1,\dots,4\}$, $q_{\max} \in \{1,\dots,5\}$, $\text{\textsc{flip}} \in \{true,false\}$, and $\text{\textsc{randroutes}} \in \{true,false\}$, the latter only when $\text{\textsc{otype}} = \mathrm{F}$. We set $maxExperiments = 1000$, with a time limit of 1200 seconds per run, and kept all remaining parameters at their default values. The calibration objective was the value of the solution found within the time limit. Within this budget, Irace evaluated 137 configurations. Because the elite configurations returned by Irace performed similarly, they were re-evaluated on the 15 calibration instances, with three runs per instance, and the one with the best average performance was adopted, namely $\text{\textsc{otype}} = \mathrm{B}$, $r_{\max} = 3$, $\text{\textsc{neighborhood\_order}} = 3$, $q_{\max} = 5$, and $\text{\textsc{flip}} = true$. Measured with respect to the best of the solutions obtained by these configurations for the same instance and seed, their average deviations ranged from 0.25\% to 0.40\%, which suggests that the performance of GVNS varies little among the best configurations found by Irace.

\subsection{Experiments with the small- and medium-scale instances of the CVRPN} \label{sec57}

We ran the GVNS metaheuristic on the 134 small- and medium-scale instances described in Section~\ref{sec52}. We adopted the parameter configuration obtained with Irace and imposed a time limit of 60 seconds per run for the 65 instances with up to 15 customers and 1200 seconds per run for the remaining ones. Because the GVNS heuristic contains randomized components, each instance was solved 10 times, the $k$-th run using the random seed $123456k$. Tables~7 and~8 in the Supplementary Material report the best result obtained over the 10 runs for each instance. In all 65 instances with up to 15 customers and in 61 of the 69 remaining ones, the best solution is found in all 10 runs, and the coefficient of variation over the 10 runs never exceeds 0.34\%. These results can be compared with the solutions obtained by CPLEX, using a time limit of one hour. As in Section~\ref{sec54}, two solution values are considered equal when their relative difference is at most $10^{-5}$. CPLEX found provably optimal solutions for 69 instances, which essentially correspond to the smallest instances in the test set. In all such cases, the GVNS heuristic also found an optimal solution. For the remaining 65 instances, CPLEX reached the time limit without proving optimality. Among these 65 instances, the GVNS heuristic matched the best solution value obtained by CPLEX in 12 cases. In the remaining 53 instances, the solution found by GVNS improved upon the best solution obtained by CPLEX. On average, these solutions were 3.16\% better. This value is the arithmetic mean, over these 53 instances, of the absolute value of the gap $100\times(\mathrm{TD}_{\textsc{gvns}} - \mathrm{TD}_{\textsc{cplex}})/\mathrm{TD}_{\textsc{cplex}}$, where TD denotes the total distance of a solution, reported in Tables~7 and~8 in the Supplementary Material. Overall, the GVNS metaheuristic matched or improved all solutions obtained by CPLEX.

\subsection{Experiments with the large-scale instances adapted from the CVRP literature} \label{sec58}

\begin{figure}[ht!]
\centering
\begin{tabular}{@{}c@{\hspace{0.09\textwidth}}c@{}}
\includegraphics[width=0.43\textwidth]{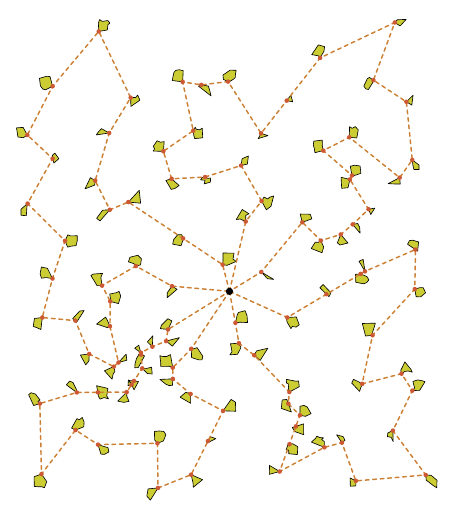} & \includegraphics[width=0.43\textwidth]{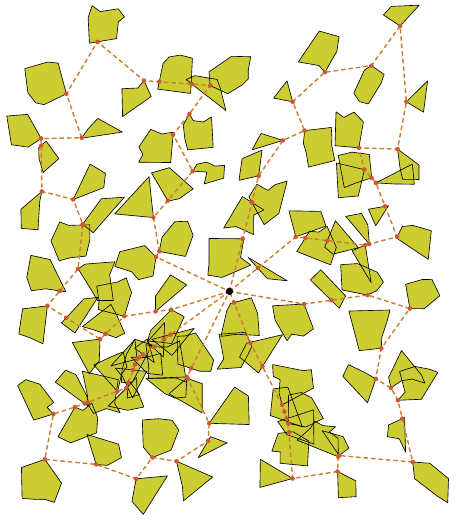} \\[0.5ex]
(a) Area of 0.1\% & (b) Area of 1.0\%
\end{tabular}

\caption{Solutions obtained for instance P-n101-k4, drawn at the same scale. Both solutions have four routes, with total distances of 611.03 and 488.65, respectively, and only 53 of the 104 edges of solution (a) are kept in solution (b).}
\label{fig5}
\end{figure}

We ran the GVNS metaheuristic on the 318 large-scale instances described in Section~\ref{sec55}. We adopted the parameter configuration obtained with Irace and imposed a time limit of 1200 seconds per run (more on this later). Because the GVNS heuristic contains randomized components, each instance was solved 10 times. Table~\ref{tabsum58} summarizes the results for each area, and Tables~9 to~14 in the Supplementary Material report the results for each instance. For each instance, we report the best total distance obtained across the 10 runs, the average total distance over all runs, and the coefficient of variation (CV), defined as the standard deviation divided by the average total distance. For 162 of the 318 instances, mostly the smaller ones, the best solution is found in all 10 runs. For the remaining 156 instances, the frequency with which the best solution is found over the 10 runs decreases as the number of customers increases. For these 156 instances, the coefficient of variation never exceeds 1.37\% and averages 0.29\%. In an ablation experiment on 30 large-scale instances (10 instances in each of the three areas, three runs of 1200 seconds each), one neighborhood was removed at a time. For every neighborhood, the average difference from the complete GVNS was below 0.2\% in absolute value, both in the best and in the mean value of the three runs, and no difference was statistically significant after the Holm correction for multiple comparisons (Table~19 in the Supplementary Material). Because removing a single neighborhood did not lead to a statistically significant improvement, the configuration was kept unchanged. As an illustration, Figure~\ref{fig5} shows the solutions obtained for instance P-n101-k4 from benchmark set P, with areas of 0.1\% and 1.0\%.

\begin{table}[ht!]
\centering
\resizebox{\textwidth}{!}{
\begin{tabular}{lrrr@{\hspace{8pt}}c@{\hspace{8pt}}rrrr@{\hspace{8pt}}c@{\hspace{8pt}}r}
\hline
 & \multicolumn{3}{c}{GVNS (10 runs)} & & \multicolumn{4}{c}{GVNS vs.\ CPLEX} & & GVNS vs.\ CVRP \\
\cline{2-4}\cline{6-9}\cline{11-11}
Area & Best in all runs & Max.\ CV (\%) & Time to best (s) & & Optimal & Ties & Better & Improv.\ (\%) & & Improv.\ (\%) \\
\hline
0.1\% & 58  & 1.06 & 284.8 & & 4 & 1 & 100 & 4.35 & & 5.33 \\
0.5\% & 52  & 1.26 & 254.9 & & 2 & 5 & 98  & 4.36 & & 10.25 \\
1.0\% & 52  & 1.37 & 321.6 & & 1 & 3 & 101 & 4.53 & & 13.34 \\
\hline
All   & 162 & 1.37 & 287.1 & & 7 & 9 & 299 & 4.41 & & 9.64 \\
\hline
\end{tabular}}
\caption{Summary of the results for the 106 large-scale instances of each area. ``Best in all runs'' counts the instances in which all 10 runs found the best solution, and ``Time to best'' is the average CPU time to find it. ``Optimal'' counts the instances solved to proven optimality by CPLEX (GVNS found all these optimal solutions), and ``Ties'' and ``Better'' count the remaining instances in which GVNS matches or improves upon CPLEX (CPLEX found no solution for Tai385). ``Improv.'' is the average improvement over CPLEX, in the instances in which GVNS is better, and over the CVRP-based upper bound.}
\label{tabsum58}
\end{table}

To validate the performance of the proposed GVNS, we solved the 318 large-scale instances using CPLEX, with a time limit of one hour; these results are reported in Tables~15 and~16 in the Supplementary Material. CPLEX achieved optimality in 7 instances and found feasible solutions, without proving optimality, for 308. In the three versions of instance Tai385, it failed to return a feasible solution due to an out-of-memory error. GVNS matched all 7 optimal solutions. Among the remaining 308 instances, GVNS tied with CPLEX in 9 and outperformed it in 299, yielding an average improvement of 4.41\% (arithmetic mean over these 299 instances). In these 299 instances, GVNS found better feasible solutions than CPLEX within the respective time limits, and the optimality of these solutions is not established.

For each CVRPN instance, we also compared the solution obtained by GVNS with the best-known solution reported in the literature for the corresponding CVRP instance from which the CVRPN instance was derived. Because the best-known solution of the original CVRP instance is also feasible for the associated CVRPN instance, it provides a valid upper bound. The results show that the solutions computed by GVNS improved upon this upper bound in all instances, with an average improvement of 9.64\%. This value is the arithmetic mean, over the 318 instances, of the absolute value of the gap $100\times(\mathrm{TD}_{\textsc{gvns}} - \mathrm{TD}^{*})/\mathrm{TD}^{*}$ reported in Tables~9 to~14 in the Supplementary Material, where $\mathrm{TD}^{*}$ denotes the total distance of the CVRP solution recomputed with unrounded Euclidean distances. As expected, the improvement grows with the size of the neighborhoods, from 5.33\% with an area of 0.1\% to 13.34\% with an area of 1.0\% (Table~\ref{tabsum58}). It should be stressed that this value does not represent an improvement of the CVRP solutions themselves: it compares objective values of two different problems and measures how much the solutions found by GVNS benefit from the possibility of choosing the visit points. When a CVRP instance is transformed into a CVRPN instance, the geographical location of each customer is replaced by a customer neighborhood containing that location. This modification enlarges the feasible region of the problem. A natural question is therefore whether the observed improvement arises merely from relocating the visit points associated with the best-known CVRP solution while preserving its routes. To answer it, we applied BCD to the routes of the best-known CVRP solution of each instance. This baseline improves upon $\mathrm{TD}^{*}$ by 4.93\%, 9.24\%, and 11.86\% on average for areas of 0.1\%, 0.5\%, and 1.0\%, which shows that BCD alone captures most of the benefit of the neighborhoods. Starting only from the solutions of the constructive heuristics, without using the CVRP solutions, GVNS matches or improves upon this baseline in 308 of the 318 instances (287 of the 288 instances with up to 100 customers), with an additional reduction, relative to the baseline, of 0.42\%, 1.13\%, and 1.71\% on average (up to 6.46\%), which grows with the size of the neighborhoods (see Figure~\ref{fig5}, in which this reduction is 1.34\% with an area of 0.1\% and 6.38\% with an area of 1.0\%). For sets A, B, E, F, M, and P, all reported CVRP solutions are proven optimal for the original CVRP instances, with rounded distances. For the Tai set, in which distances are not rounded, all solutions are proven optimal except for instance Tai385.

\begin{figure}[ht!]
\centering
\begingroup
  \makeatletter
  \providecommand\color[2][]{%
    \GenericError{(gnuplot) \space\space\space\@spaces}{%
      Package color not loaded in conjunction with
      terminal option `colourtext'%
    }{See the gnuplot documentation for explanation.%
    }{Either use 'blacktext' in gnuplot or load the package
      color.sty in LaTeX.}%
    \renewcommand\color[2][]{}%
  }%
  \providecommand\includegraphics[2][]{%
    \GenericError{(gnuplot) \space\space\space\@spaces}{%
      Package graphicx or graphics not loaded%
    }{See the gnuplot documentation for explanation.%
    }{The gnuplot epslatex terminal needs graphicx.sty or graphics.sty.}%
    \renewcommand\includegraphics[2][]{}%
  }%
  \providecommand\rotatebox[2]{#2}%
  \@ifundefined{ifGPcolor}{%
    \newif\ifGPcolor
    \GPcolortrue
  }{}%
  \@ifundefined{ifGPblacktext}{%
    \newif\ifGPblacktext
    \GPblacktexttrue
  }{}%
  \let\gplgaddtomacro\g@addto@macro
  \gdef\gplbacktext{}%
  \gdef\gplfronttext{}%
  \makeatother
  \ifGPblacktext
    \def\colorrgb#1{}%
    \def\colorgray#1{}%
  \else
    \ifGPcolor
      \def\colorrgb#1{\color[rgb]{#1}}%
      \def\colorgray#1{\color[gray]{#1}}%
      \expandafter\def\csname LTw\endcsname{\color{white}}%
      \expandafter\def\csname LTb\endcsname{\color{black}}%
      \expandafter\def\csname LTa\endcsname{\color{black}}%
      \expandafter\def\csname LT0\endcsname{\color[rgb]{1,0,0}}%
      \expandafter\def\csname LT1\endcsname{\color[rgb]{0,1,0}}%
      \expandafter\def\csname LT2\endcsname{\color[rgb]{0,0,1}}%
      \expandafter\def\csname LT3\endcsname{\color[rgb]{1,0,1}}%
      \expandafter\def\csname LT4\endcsname{\color[rgb]{0,1,1}}%
      \expandafter\def\csname LT5\endcsname{\color[rgb]{1,1,0}}%
      \expandafter\def\csname LT6\endcsname{\color[rgb]{0,0,0}}%
      \expandafter\def\csname LT7\endcsname{\color[rgb]{1,0.3,0}}%
      \expandafter\def\csname LT8\endcsname{\color[rgb]{0.5,0.5,0.5}}%
    \else
      \def\colorrgb#1{\color{black}}%
      \def\colorgray#1{\color[gray]{#1}}%
      \expandafter\def\csname LTw\endcsname{\color{white}}%
      \expandafter\def\csname LTb\endcsname{\color{black}}%
      \expandafter\def\csname LTa\endcsname{\color{black}}%
      \expandafter\def\csname LT0\endcsname{\color{black}}%
      \expandafter\def\csname LT1\endcsname{\color{black}}%
      \expandafter\def\csname LT2\endcsname{\color{black}}%
      \expandafter\def\csname LT3\endcsname{\color{black}}%
      \expandafter\def\csname LT4\endcsname{\color{black}}%
      \expandafter\def\csname LT5\endcsname{\color{black}}%
      \expandafter\def\csname LT6\endcsname{\color{black}}%
      \expandafter\def\csname LT7\endcsname{\color{black}}%
      \expandafter\def\csname LT8\endcsname{\color{black}}%
    \fi
  \fi
    \setlength{\unitlength}{0.0500bp}%
    \ifx\gptboxheight\undefined%
      \newlength{\gptboxheight}%
      \newlength{\gptboxwidth}%
      \newsavebox{\gptboxtext}%
    \fi%
    \setlength{\fboxrule}{0.5pt}%
    \setlength{\fboxsep}{1pt}%
    \definecolor{tbcol}{rgb}{1,1,1}%
\begin{picture}(9070.00,4648.00)%
    \gplgaddtomacro\gplbacktext{%
      \csname LTb\endcsname
      \put(820,1653){\makebox(0,0)[r]{\strut{}1.00}}%
      \put(820,2156){\makebox(0,0)[r]{\strut{}1.01}}%
      \put(820,2660){\makebox(0,0)[r]{\strut{}1.02}}%
      \put(820,3163){\makebox(0,0)[r]{\strut{}1.03}}%
      \put(820,3666){\makebox(0,0)[r]{\strut{}1.04}}%
      \put(820,4169){\makebox(0,0)[r]{\strut{}1.05}}%
      \put(1029,1232){\makebox(0,0){\strut{}$0$}}%
      \put(1543,1232){\makebox(0,0){\strut{}$10$}}%
      \put(2057,1232){\makebox(0,0){\strut{}$20$}}%
      \put(2571,1232){\makebox(0,0){\strut{}$30$}}%
      \put(3084,1232){\makebox(0,0){\strut{}$40$}}%
      \put(3598,1232){\makebox(0,0){\strut{}$50$}}%
      \put(4112,1232){\makebox(0,0){\strut{}$60$}}%
    }%
    \gplgaddtomacro\gplfronttext{%
      \csname LTb\endcsname
      \put(162,2961){\rotatebox{-270.00}{\makebox(0,0){\strut{}Average of $\mathrm{TD}(t)/\mathrm{TD}_{\mathrm{best}}$}}}%
      \put(2570,572){\makebox(0,0){\strut{}\begin{tabular}{c}CPU time (min)\\[1.6ex](a) The 288 instances with up to\\100 customers (one-hour runs)\end{tabular}}}%
    }%
    \gplgaddtomacro\gplbacktext{%
      \csname LTb\endcsname
      \put(5518,1653){\makebox(0,0)[r]{\strut{}1.00}}%
      \put(5518,2156){\makebox(0,0)[r]{\strut{}1.01}}%
      \put(5518,2660){\makebox(0,0)[r]{\strut{}1.02}}%
      \put(5518,3163){\makebox(0,0)[r]{\strut{}1.03}}%
      \put(5518,3666){\makebox(0,0)[r]{\strut{}1.04}}%
      \put(5518,4169){\makebox(0,0)[r]{\strut{}1.05}}%
      \put(5728,1232){\makebox(0,0){\strut{}$0$}}%
      \put(6353,1232){\makebox(0,0){\strut{}$60$}}%
      \put(6978,1232){\makebox(0,0){\strut{}$120$}}%
      \put(7604,1232){\makebox(0,0){\strut{}$180$}}%
      \put(8229,1232){\makebox(0,0){\strut{}$240$}}%
      \put(8854,1232){\makebox(0,0){\strut{}$300$}}%
    }%
    \gplgaddtomacro\gplfronttext{%
      \csname LTb\endcsname
      \put(4860,2961){\rotatebox{-270.00}{\makebox(0,0){\strut{}Average of $\mathrm{TD}(t)/\mathrm{TD}_{\mathrm{best}}$}}}%
      \put(7291,572){\makebox(0,0){\strut{}\begin{tabular}{c}CPU time (min)\\[1.6ex](b) The 30 instances with more than\\100 customers (five-hour runs)\end{tabular}}}%
    }%
    \gplbacktext
    \put(0,0){\includegraphics[width={453.50bp},height={232.40bp}]{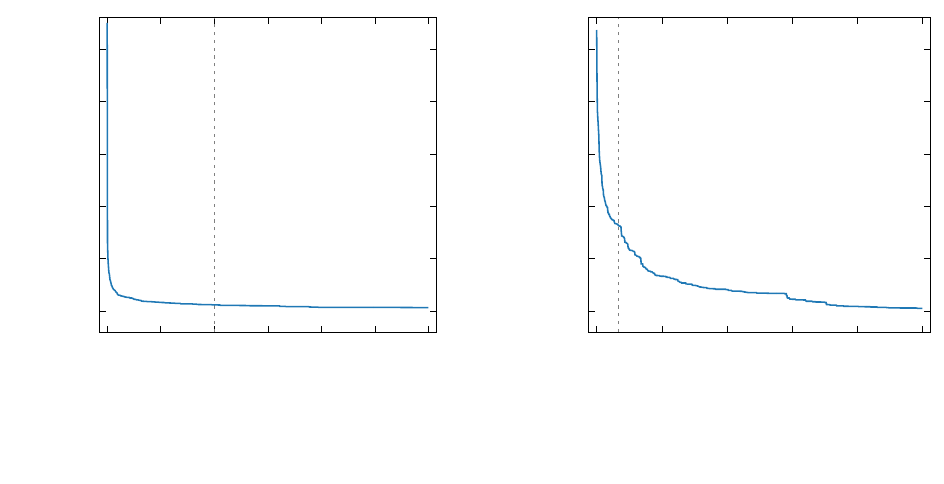}}%
    \gplfronttext
  \end{picture}%
\endgroup

\caption{Average of $\mathrm{TD}(t)/\mathrm{TD}_{\mathrm{best}}$, where $\mathrm{TD}(t)$ is the total distance of the incumbent solution at time $t$ and $\mathrm{TD}_{\mathrm{best}}$ is the best total distance found for the instance. The dashed line marks 20 minutes (1200 seconds).}
\label{fig6}
\end{figure}

The 1200-second time limit adopted in the computational experiments of this section was chosen arbitrarily. Because it is unclear whether GVNS could find better-quality solutions with additional computational time, we reran GVNS on all 318 instances, with a time limit of one hour for the 288 instances with up to 100 customers and five hours for the 30 instances with more than 100 customers. Because this experiment was computationally expensive, each instance was solved only once. Tables~17 and~18 in the Supplementary Material summarize the computational results. For the instances with up to 100 customers, the solution improves after 1200 seconds in only 32 of the 288 instances, with an average improvement of 0.48\% among them. In contrast, for the instances with more than 100 customers, better solutions are found in 27 of the 30 instances, with an average improvement of 1.71\%.

Figure~\ref{fig6} shows the evolution of the average incumbent solution value as a function of time for the instances with up to 100 customers and for those with more than 100 customers. To prevent instances with larger objective values from exerting a disproportionate influence on the averages, the incumbent value of each instance was normalized by dividing it by the best solution found for that instance, either in the long run or in the 10 runs of 1200 seconds. Figure~\ref{fig6}(a) indicates that the 1200-second time limit was sufficient for the instances with up to 100 customers. Improvements after the first 1200 seconds were small, and the average incumbent curve quickly flattens, suggesting that near-convergence is reached within the allotted time. In contrast, Figure~\ref{fig6}(b) indicates that a 1200-second time limit is insufficient for the larger instances, as the average incumbent curve continues to decrease throughout the five-hour horizon. This behavior suggests that the search remains capable of finding improving solutions even after several hours of computation and that additional gains could potentially be achieved through longer runs or stronger intensification strategies. Consequently, for the instances with more than 100 customers, the results obtained with the 1200-second limit should be read as what GVNS achieves within that time, and not as the best solutions it can find; this is a limitation of the adopted experimental protocol.

In the routes of the initial solutions of the large-scale instances, BCD starting from the barycenters required on average about 11, 19, and 26 cycles per route for areas of 0.1\%, 0.5\%, and 1.0\%, respectively, and never reached the maximum of 1000 cycles. Finally, inside the search, BCD starts from the visit points of the current solution. In the 10 runs of each instance in Sections~\ref{sec57} and~\ref{sec58}, it required on average 2.3 and 11.2 cycles per call, respectively (averages over the instances), and reached the maximum of 1000 cycles in fewer than 0.002\% of the calls.

\section{Conclusions and future work} \label{sec6}

In this work, we studied a capacitated vehicle routing problem with polygonal neighborhoods (CVRPN) in which customers are represented by simple polygons, possibly non-convex. In contrast with classical VRP formulations, a route is defined not only by the sequence of customers to be visited but also by the selection of a visit point within each customer neighborhood. This additional geometric component substantially increases the complexity of the problem and requires the repeated solution of an MIQCP during the search process.

To address this challenge, we proposed a block coordinate descent (BCD) procedure for computing visit points associated with a given visit sequence. Computational experiments showed that, in the instances tested, BCD provides solutions of the same quality as those obtained by solving the MIQCP formulations with a commercial solver while requiring only a tiny fraction of the computational effort, often reducing the solution time by several orders of magnitude. This result makes the repeated evaluation of candidate routes feasible within a metaheuristic framework. Therefore, a central contribution of the proposed framework is the BCD route-evaluation procedure, which enables the efficient optimization of visit points during the search process and makes large-scale GVNS exploration computationally tractable.

Building on this route-evaluation mechanism, we developed a GVNS metaheuristic tailored to the CVRPN. The method combines several intra-route and inter-route search neighborhood structures and was calibrated automatically using Irace. Extensive experiments on both randomly generated instances and large-scale instances derived from CVRPLib benchmarks demonstrate the effectiveness of the proposed approach. For the small- and medium-scale instances, the method matched all solutions proven optimal by CPLEX and frequently improved the best solutions obtained by the exact solver within the prescribed time limit. For the large-scale benchmark set, the proposed method improved the CVRP-based upper bound in all instances, yielding an average improvement of 9.64\%, which grows with the size of the neighborhoods. Most of this improvement is already obtained by applying BCD to the best-known CVRP routes, and GVNS, starting only from constructive heuristics, matches or improves upon this baseline in 308 of the 318 instances.

The analysis of the incumbent trajectories of the large-scale instances reveals a marked difference depending on the number of customers. For the instances with up to 100 customers, the average incumbent value stabilizes within the 1200-second (20-minute) horizon considered in the experiments, whereas for those with more than 100 customers it continues to decrease throughout the entire 18,000-second (5-hour) horizon. These results suggest that the proposed GVNS approaches convergence on the instances with up to 100 customers, whereas additional computational effort or stronger intensification strategies may still yield further improvements on the larger ones.

Overall, the results indicate that the combination of GVNS and BCD constitutes an effective framework for solving routing problems with geometric customer neighborhoods. Moreover, the proposed BCD procedure is independent of the surrounding metaheuristic and may therefore be incorporated into alternative search strategies. The mathematical model, benchmark instances, and software developed in this work provide a foundation for future research on routing problems with customer neighborhoods.

The proposed approach has limitations that should be acknowledged. BCD does not guarantee global optimality of the visit points and, in the simplified version adopted here, it may cycle in theory, although this was not observed in our experiments. In addition, because $\Psi$ is not differentiable where consecutive visit points coincide, which may happen when neighborhoods overlap, BCD may stop at points that are not local minimizers of~\eqref{probb}. The largest instances remain computationally challenging: for the instances with more than 100 customers, the average incumbent value was still decreasing after five hours of computation. Moreover, the customer neighborhoods of all the instances considered were generated artificially, either at random or around the customers of CVRP instances from the literature, so the performance of the method on real-world geographical data remains to be established. Also, the three CVRPN instances derived from each CVRP instance have neighborhoods of different sizes, and a single realization of the neighborhoods was generated for each size, so the results may depend on the particular generation process.

Several directions deserve further investigation. Future work includes the development of exact and matheuristic approaches, the study of richer variants involving time windows, heterogeneous fleets, or multiple depots, and the extension of the framework to three-dimensional instances motivated by aerial and robotic applications. Another promising direction is the design of specialized search neighborhood structures that jointly modify customer sequences and visit-point locations. The experiments with three neighborhood sizes show that the savings with respect to the classical CVRP grow with the size of the neighborhoods. It would also be interesting to investigate systematically how this size, relative to the distances between customers, affects the difficulty of the instances. Finally, the incumbent trajectories observed for large-scale instances suggest that investigating longer-run search strategies, adaptive stopping criteria, and more powerful intensification mechanisms may lead to further improvements in solution quality.\\

\noindent
\textbf{Acknowledgments:} We sincerely thank Marcelo Marazzi for his careful editing and helpful suggestions that enhanced the quality of the written text in English. This work was partially supported by the Brazilian agencies FAPESP [grant numbers \mbox{2022/05803-3} and \mbox{2023/08706-1}] and CNPq [grant number \mbox{302073/2022-1}].

\medskip
\noindent
\textbf{Declaration of generative AI use:} During the preparation of this manuscript, the authors used Claude Opus 5.5 (Anthropic) for final language editing. The authors reviewed all the suggested changes and take full responsibility for the content of this publication.

\printbibliography

\end{document}